\documentclass[11pt,a4paper]{article}
\usepackage[utf8]{inputenc}
\usepackage[margin=2cm]{geometry}
\usepackage{bm} 
\usepackage{enumitem} 
\usepackage{array}  
\usepackage[pdftex]{graphicx} 
\usepackage{subcaption} 
\usepackage{amssymb} 
\usepackage{amsmath} 
\usepackage{amsthm} 
\usepackage{natbib}
\usepackage{hyperref} 
\usepackage{authblk} 
\usepackage{fancyhdr} 
\usepackage{mathtools} 
\usepackage{tikz-cd}
\usepackage{multirow}
\usepackage{lineno}

\newcommand{\pfrac}[2]{\frac{\partial{#1}}{\partial{#2}}}
\newcommand{\dx}[1]{\hspace{0.75mm}\mathrm{d}{#1}}

\fancypagestyle{plain}{
\fancyhf{}
\lfoot{\textcopyright \ Crown Copyright, Met Office}
\cfoot{\thepage}
}

\providecommand{\keywords}[1]
{
  \small	
  \textbf{\textit{Keywords---}} #1
}

\title{SWIFT: A Monotonic, Flux-Form Semi-Lagrangian Tracer Transport Scheme for Flow with Large Courant Numbers}

\author[1]{Thomas M. Bendall\thanks{Corresponding author, contact at thomas.bendall@metoffice.gov.uk}}
\author[1]{James Kent}
\affil[1]{Dynamics Research, Met Office, Exeter, UK}
\date{}

\begin{document}

\maketitle

\begin{abstract}
\noindent
Local conservation of mass and entropy are becoming increasingly desirable properties for modern numerical weather and climate models.
This work presents a Flux-Form Semi-Lagrangian (FFSL) transport scheme, called SWIFT, that facilitates this conservation for tracer variables, whilst maintaining other vital properties such as preservation of a constant, monotonicity and positivity.
Importantly, these properties all hold for large Courant numbers and multi-dimensional flow, making the scheme appropriate for use within a dynamical core which takes large time steps. \\
\\
The SWIFT scheme presented here can be seen as an evolution of the FFSL methods of Leonard et al and Lin and Rood.
Two-dimensional and three-dimensional schemes consist of a splitting into a sequence of one-dimensional calculations.
The new SWIFT splitting presented here allows monotonic and positivity properties from the one-dimensional calculations to be inherited by the multi-dimensional scheme.
These one-dimensional calculations involve separating the mass flux into terms that correspond to integer and fractional parts of the Courant number.
Key to achieving conservation is coupling the transport of tracers to the transport of the fluid density, through re-use of the discrete mass flux that was calculated from the fluid density in the transport of the tracers.
This work also describes how these properties can still be attained when the tracer is vertically-staggered from the density in a Charney-Phillips grid.

\end{abstract}

\section{Motivation and Overview} \label{sec:intro}
\subsection{Desirable Properties for Transport Schemes}
Numerical methods for solving the transport equation are a core component of atmospheric models that are used for weather prediction and climate research.
These equations typically take two forms:
\begin{subequations}
\begin{align}
& \pfrac{\rho}{t} + \bm{\nabla\cdot}\left(\rho\bm{u}\right) = 0, \label{eqn:conservative_transport} \\
& \pfrac{m}{t} + \left(\bm{u\cdot\nabla}\right)m = 0, \label{eqn:advective_transport}
\end{align}
\end{subequations}
which are known as the \textit{conservative} form \eqref{eqn:conservative_transport} for a density field $\rho$ or the \textit{advective} form \eqref{eqn:advective_transport} which may be used for instance for a tracer mixing ratio\footnote{The mixing ratio of a species is defined as the ratio of the mass density of the species to some reference mass density, which in this work is denoted by $\rho$. In atmospheric applications this could be the density of dry air.} $m$. 
In both forms, $\bm{u}$ is the wind velocity field.
Considering the equation for the product $\rho m$ (which describes the density of the tracer species), it can be shown that this obeys:
\begin{equation}
\pfrac{\rho m}{t} + \bm{\nabla\cdot}\left(\rho m\bm{u}\right) = 0. \label{eqn:tracer_conservative}
\end{equation}
The relevance of this form of the transport equation is that discretisations of \eqref{eqn:conservative_transport} and \eqref{eqn:tracer_conservative} provide local (and hence global) mass conservation, in that each cell has a closed mass budget when taking into account the mass fluxes into neighbouring cells.
As discussed by \cite{thuburn2008some}, these conservation properties are desirable for weather and particularly climate models, so that the discretisation itself does not become a net source or sink of the tracer.
Locally-conservative schemes naturally provide a closed budget for the tracer species. 
Additionally, for the transport of a density, \citet{lin1996ffsl} and \citet{melvin2024mixed} argue that any conservative transport scheme should provide increments that are linear in $\bm{\nabla\cdot u}$ when the density is constant, a property that we will describe here as \textit{divergence-linearity}.
The corollary of this is that any density field remains constant under divergence-free flow.\footnote{This property is described as ``constancy'' by \citet{leonard1996cosmic}, while it is a sub-property of the ``consistency'' described by \citet{lin1996ffsl}.}
\\
\\
However other properties are also important for tracer transport schemes: it can be vital that the scheme is \textit{positivity-preserving} so that unphysical negative concentrations are not produced.
A stronger possible requirement is that of \textit{monotonicity}, which prevents the introduction of new minima or maxima in the solution.
Note that as discussed by \citet{lauritzen2011numerical}, in contrast to $m$, the product $\rho m$ will not in general remain bounded by its initial extrema in divergent flows, but its evolution should be \textit{shape-preserving} in the sense that no new unphysical extrema are generated.
In this work we focus on the monotonicity requirement of keeping the minima and maxima of mixing ratios bounded by their initial values.
It is also important for a conservative tracer transport scheme to be consistent with both the advective form of the equation and the conservative transport of $\rho$ (see \citep{jockel2001fundamental} and \citep{zhang2008consistency}).
Thus the property of \textit{consistency} is defined as the preservation of a constant mixing ratio, independent of the flow.
Finally, many operational numerical weather prediction models use large time steps.
It is therefore highly attractive to use a scheme which is stable at large Courant numbers\footnote{For instance, the finest resolution global configurations discussed in \citet{walters2019met} correspond to Courant numbers of around 5 in the mid-latitude jets.}
whilst maintaining computational efficiency and without compromising the other properties listed here.
A complete discussion of properties that are desirable for transport schemes can be found in \citet{lauritzen2011numerical}.
\\
\\
In addition to the desirable properties listed above, the transport scheme must work within the parameters of the atmospheric dynamical core, which may result in additional requirements for the transport scheme.
Many dynamical cores (e.g. \citep{chen2008new,wood2014inherently,girard2014staggered,melvin2024mixed}) use hexahedral cells and make use of C-grid staggering in the horizontal and Charney-Phillips staggering \citep{charney1953verticalstaggering} in the vertical, primarily due to the favourable representation of gravity wave dispersion (see for instance \citet{staniforth2012horizontal} and \citet{thuburn2022numerical}). 
With these staggerings, the density $\rho$ of dry air is represented at cell centres, the wind velocity $\bm{u}$ is described at the centres of cell faces and the potential temperature $\theta$ is co-located with the vertical component of the wind velocity (at the centre of the top/bottom faces of cells).
It is common for tracers such as moisture mixing ratios to be co-located with the potential temperature (or other entropy-type variable) as this is beneficial for capturing the latent heat exchanges associated with phase changes \citep{bendall2023solution}.
This work considers a transport scheme that can be used with these staggerings.

\subsection{Background to Flux-Form Semi-Lagrangian Schemes}
Some notable conservative transport methods which meet many of the criteria above include: the CSLAM scheme of \citet{lauritzen2010conservative} and \citet{harris2011flux}; the SLICE scheme of \citet{zerroukat2002slice,zerroukat2004slice} and \citet{zerroukat2012three}; the multi-moment scheme of \cite{tang2022three}; and the substepped finite volume Method-of-Lines method described in \citet{melvin2024mixed}.
However, the focus in this work is another class of schemes that possesses many of the desirable properties: Flux-Form Semi-Lagrangian (FFSL) schemes.
These schemes originate from the flux integration methods of \citet{van1974towards,colella1984ppm,carpenter1990application}
and were developed in multiple dimensions by \citet{leonard1996cosmic} and \citet{lin1996ffsl}.
Thanks to their conservation and stability properties, FFSL schemes have been successfully used for transport in geophysical fluid models \citep{lin2004fvcore,putman2007fvtransport,zhou2012computational,neale2013mean,gillibrand2016mass,harris2021scientific,zhang2023history,mouallem2023implementation}.
FFSL schemes involve computing mass fluxes by integrating fields between flux points and corresponding departure points.
Conservation is inherent, as increments are evaluated from the divergence of fluxes, which are shared between cells.
Like standard semi-Lagrangian schemes, stability is not conditional on the Courant number which permits large time steps. \\
\\
Instead, the stability of the scheme depends upon the Lipschitz number, which measures the one-dimensional divergence of the velocity $\bm{u}$ seen by a discretisation with time step $\Delta t$.
In other words $\Delta t|\partial u /\partial x|$ in the $x$-direction and similar for other directions.
In an Arakawa C-grid staggering, in which density fields are located at cell centres and the velocity is described by its components normal to each facet, the Lipschitz numbers are naturally calculated at cell centres.
However Lipschitz numbers can also be calculated at facets, so that in one-dimension $\ell=(c-c^\dagger) \, \mathrm{sign}(c)$, where $c$ is the Courant number at that facet and $c^\dagger$ is the Courant number at the next facet upwind.
The result is that when this number is above 1, the trajectories of discrete particles can cross \citep{smolarkiewicz1992class}.\\
\\
The multi-dimensional forms of FFSL presented by \citet{leonard1996cosmic} and \citet{lin1996ffsl} are built from one-dimensional calculations through dimensional-splitting techniques.
Following \citet{leonard1996cosmic}, the multi-dimensional splitting of \citet{leonard1996cosmic} and \citet{lin1996ffsl} will be referred to as the COSMIC splitting throughout this work.
Advantages of using multiple one-dimensional methods over pure multi-dimensional methods include the ease of applying limiters, the ease of extending to large Courant number flow, and that in general one-dimensional methods are computationally cheaper than their multi-dimensional counterparts. \\
\\
However one drawback of the COSMIC scheme, as noted by both \citet{leonard1996cosmic} and \citet{lin1996ffsl},  is that applying monotonic limiters to the one-dimensional calculations does not guarantee strict monotonicity in the multi-dimensional scheme.
For flows with Courant numbers less than 1 and with the use of one-dimensional limiters in the COSMIC scheme, the minima and maxima of fields rarely exceed the initial values.
Thus COSMIC was called ``essentially shape-preserving'' by \citet{leonard1996cosmic} and ``quasi-monotonic'' by \citet{putman2007fvtransport}.
However at large Courant numbers, COSMIC does not have monotonicity \citep{leonard1996cosmic,bott2010improving}.
Multi-dimensional limiters do exist that do not cause loss of conservation (for instance that of \citet{thuburn1996multidimensional}) but as of writing, the authors do not know of any that are strictly monotonic for COSMIC splitting with Courant numbers greater than 1 without compromising another desirable property.
Therefore if strict monotonic or positive-definite properties are required, the COSMIC scheme is generally substepped, losing the benefits of the stability of the FFSL approach, or clipped, losing conservation. \\
\\
A final important work to note is that of \citet{skamarock2006limiters}, which applied monotonic and positivity-enforcing filters to FFSL schemes for tracer transport, with flows of large Courant numbers.
\citet{skamarock2006limiters} obtained these properties by considering simple first-order splittings in time (which have reduced accuracy in time), and also Strang splittings.
When applied to three-dimensional transport, neither of these schemes has symmetry between the two horizontal dimensions.
 
\subsection{Overview}
This work presents a transport scheme, SWIFT (Splitting With Improved FFSL for Tracers), that satisfies all of the properties discussed above and can be used for transporting density and tracer mixing ratio fields.
SWIFT is a dimensionally-split FFSL scheme that can be seen as an evolution of that of \citet{leonard1996cosmic} and \citet{lin1996ffsl}, using multiple one-dimensional steps to solve the multi-dimensional transport equation.
Similar to the approach of \cite{skamarock2006limiters}, fluxes are explicitly separated into parts corresponding to integer and fractional Courant numbers.
The subgrid reconstruction that is used for the fractional part is recast as a linear sum in the field values.
This is crucial to obtaining the consistency property, and can also be exploited to adapt low-order correction schemes to enforce positivity or monotonicity for large Courant number flows. 
However, a different two-dimensional splitting is introduced to the splitting used by the \citet{leonard1996cosmic} and \citet{lin1996ffsl} schemes.
While the multi-dimensional splitting of \citet{leonard1996cosmic} and \citet{lin1996ffsl} does not enforce strict monotonicity or positivity, the alternative splitting used by SWIFT does enforce monotonicity and positivity in multi-dimensions (provided monotonic and positivity filters are used by the one-dimensional scheme) whilst maintaining the other properties.
As well as one-dimensional meshes, the scheme is designed for use with non-orthogonal two- and three-dimensional meshes that respectively use quadrilateral and hexahedral cells.
\\
\\
The main novel aspects of this work can be summarised as:
\begin{enumerate}
\item an evolution of the work presented by \citet{skamarock2006limiters} for tracer transport, through recomputation of departure points and the use of COSMIC splitting;
\item a novel two-dimensional splitting, which preserves the monotonicity or positivity properties of the one-dimensional steps for tracers. This splitting can be applied to density variables obeying the conservative transport equation, tracers obeying a non-conservative transport equation, and tracers obeying a conservative transport equation which is consistent with the transport of some underlying density;
\item a three-dimensional version of the scheme, using Strang splitting between the one-dimensional vertical and two-dimensional horizontal directions;
\item the adaptation of these schemes to the Charney-Phillips grid, whilst maintaining the other desirable properties.
\end{enumerate}
The remainder of this work is laid out as follows.
Section \ref{sec:ffsl} reviews the one-dimensional flux calculations used by FFSL schemes, and the two-dimensional  splitting used by \citet{leonard1996cosmic} and \citet{lin1996ffsl}, and demonstrates how these methods can be used for conservative tracer transport. The new monotonicity-preserving splitting used by SWIFT is presented in Section \ref{sec:swift_splitting}, which also includes analysis that demonstrates that this splitting possesses the desirable properties described in this section.
Results from testing the splitting in two-dimensions are presented in Section \ref{sec:results2d}.
A three-dimensional form of this scheme is presented in Section \ref{sec:3D_CP}, which uses a Strang splitting for the vertical and horizontal steps. This section also shows how this scheme is adapted to the Charney-Phillips grid, using the framework of \cite{bendall2023solution}.

\section{Flux-Form Semi-Lagrangian Transport} \label{sec:ffsl}

\subsection{Review of One Dimensional Flux-Form Semi-Lagrangian Schemes} \label{sec:1d_ffsl}

\citet{lin1996ffsl} describe large time-step-permitting finite-volume schemes as one-dimensional FFSL schemes. These one-dimensional FFSL schemes compute the flux by integrating the field from the flux point to a departure point. The Courant number is separated into integer and fractional parts. The integer part corresponds to the integral of the field over entire grid cells, which is just the sum over these cells of the field multiplied by the cell volume. The fractional part corresponds to the spatial integral of the field within the cell in which the departure point lies, to compute the mass corresponding to the fractional part of the Courant number. A subgrid reconstruction is used to represent the field within the cell, and this is integrated from the cell edge to the departure point. \\
\\
In this section, the one-dimensional FFSL calculation is reviewed, but on a general one-, two- or three-dimensional mesh.
At first, calculations of mass fluxes are described for a density $\rho$.
For simplicity, the one-dimensional scheme is described in the $x$ direction only, where $x_i$ are the coordinates of cell centres, $i$ is the grid index, $V_i$ is the volume of the $i$-th cell, and $\rho_i=\rho(x_i)$ is the density in cell $i$.
Mass fluxes are computed at cell facets, and denoted by $F_{i+1/2}$.
The cell facets $\varGamma_{i+1/2}$ are vertices for one-dimensional meshes, edges for two-dimensional meshes, and faces for three-dimensional meshes, and are indexed by an additional $1/2$ to indicate that these are staggered from cell centres.
While $\bm{u}$ is the general wind velocity vector, $u_{i+1/2}$ is its projection onto the normal $\bm{n}$ to the facet $\varGamma_{i+1/2}$.
The facet has length or area $S_{i+1/2}$, so that
\begin{equation} \label{eqn:u_normal}
u_{i+1/2} S_{i+1/2} =
\int_{\varGamma_{i+1/2}} \bm{u \cdot n} \ \mathrm{d}S = \int_{\varGamma_{i+1/2}} \bm{u} \bm{\cdot} \mathrm{d}\bm{S},
\end{equation}
where $\mathrm{d}\bm{S}$ is the measure denoting integration over the facet.
When the mesh is one-dimensional, $V_i=\Delta x_i$ and $S_{i+1/2}=1$.
For a multi-dimensional mesh, $\bm{u}$ can naturally be expanded through either primal basis vectors (which point tangential to mesh lines) or dual basis vectors (which are normal to mesh lines).
These are not equivalent if the mesh is not orthogonal, in which case the quantity $u_{i+1/2}$ corresponds to the contravariant component of $\bm{u}$, which multiplies the primal basis vector on $\varGamma_{i+1/2}$.

\subsubsection{Subgrid Reconstruction}

The subgrid reconstruction determines the accuracy of the scheme, both in terms of formal order-of-accuracy and in terms of monotonicity or positivity. A constant reconstruction results in a first-order scheme, whereas a linear reconstruction, such as that of \citet{fromm1968method} and \citet{van1974towards}, is second-order. The formally third-order quadratic reconstruction, such as used in the NIRVANA scheme of \citet{leonard1995nirvana} or the Piecewise Parabolic Method of \citet{colella1984ppm}, PPM (see Appendix A for the discretization used in this paper),  is of the form
\begin{equation} \label{eqn:quad_coef}
    \mathcal{Q}(\xi) = a_0 + a_1 \xi + a_2 \xi^2, \quad \xi \in [0, 1],
\end{equation}
where the coefficients $a_i$ are determined by the choice of scheme and are functions of the field at neighbouring grid cells. $\xi$ is the coordinate of the reference cell in the chosen direction. 

\subsubsection{Departure Point Calculation} \label{sec:dep_pt}
\begin{figure}[h!]
\centering
\includegraphics[width=0.5\textwidth]{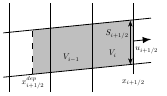}
\caption{An illustration on a two-dimensional mesh of the process used to compute the coordinate of the departure point $x^{dep}_{i+1/2}$, corresponding to the flux point $x_{i+1/2}$. The process involves finding the volume that is swept along coordinate lines over time interval $\Delta t$ through the facet of area $S_{i+1/2}$ at $x_{i+1/2}$.
This swept volume is shaded in grey.
The departure point is such that the swept volume is equal to $u_{i+1/2}S_{i+1/2}\Delta t$, thus ensuring that \eqref{eqn:dep_point} is satisfied.}
\label{fig:dep_point}
\end{figure}

The flux is computed as the integral of the mass from the flux point $x_{i+1/2}$ to the departure point $x^{dep}_{i+1/2}$.
The departure point therefore needs determining for each facet in order to compute the flux.
We use an Eulerian departure point calculation, which as will be seen ensures that if a density field is constant, any increments are proportional to $\bm{\nabla\cdot u}$.
In this case, the departure point $x^{dep}_{i+1/2}$ corresponding to facet $\varGamma_{i+1/2}$ is determined by the relationship
\begin{equation} \label{eqn:dep_point}
\int_{\varGamma_{i+1/2}} \bm{u \cdot} \mathrm{d}\bm{S} = \frac{1}{\Delta t} \int_{V(x=x^{dep}_{i+1/2})}^{V(x=x_{i+1/2})} \dx{V},
\end{equation}
where the volume integral on the right-hand side is computed along the coordinate lines used to define the mesh.
This process is illustrated in Figure \ref{fig:dep_point}. \\
 \\
The departure displacement $\delta_{i+1/2}$ is denoted by
\begin{equation}
\delta_{i+1/2} = x_{i+1/2} - x^{dep}_{i+1/2}.
\end{equation}
It can also be convenient to cast the departure displacement in a dimensionless form $c_{i+1/2}$, which is equivalent to the number of cells over which to integrate the mass field. \\
\\
For a uniform grid with grid spacing $\Delta x$ in the $x$ direction, the dimensionless departure displacement from the flux point $i+1/2$ is
\begin{equation}
c_{i+1/2} = u_{i+1/2} \Delta{t} / \Delta{x}.
\end{equation}
This demonstrates that $c_{i+1/2}$ is just the one-dimensional Courant number evaluated at the flux point.
For the flux integral calculation, the integer part is the number of whole cells to sum the field over, and the fractional part is the fraction of the departure cell to integrate the reconstruction. The Courant number can be written as
\begin{equation}
c_{i+1/2} = c^{int}_{i+1/2} + c^{frac}_{i+1/2},
\end{equation}
where $c^{int}$ is the integer part, and $c^{frac}$ the fractional part. \\
\\
For a general non-uniform grid, $c_{i+1/2}$ must be computed so that the integer and fractional parts have the same meaning as for the uniform grid.
The integer part is thus computed as
\begin{equation} \label{eqn:int_courant}
c^{int}_{i+1/2} = 
\begin{cases}
\max\left(N\in\mathbb{Z} : |u_{i+1/2}| S_{i+1/2}\Delta t \geq\sum_{j=1}^N V_{i+1-j} \right), & u_{i+1/2} \geq 0, \\
-\max\left(N\in\mathbb{Z} : |u_{i+1/2}| S_{i+1/2} \Delta t \geq\sum_{j=1}^N V_{i+j} \right), & u_{i+1/2} < 0.
\end{cases}
\end{equation}
The velocity can also be decomposed into integer and fractional parts, although only the fractional part will be used later.
This is expressed as $\widehat{u}^{frac}$, where the hat adornment $\widehat{\cdot}$ denotes that the fractional part is computed relative to the Courant number decomposition.
The fractional velocity is found from
\begin{equation} \label{eqn:frac_wind}
\widehat{u}^{frac}_{i+1/2} = 
\begin{dcases}
u_{i+1/2}-\sum_{j=1}^{c^{int}_{i+1/2}} V_{i+1-j} / \left(S_{i+1/2}\Delta t\right), & u_{i+1/2} \geq 0, \\
u_{i+1/2}+\sum_{j=1}^{\left|c^{int}_{i+1/2}\right|} V_{i+j} / \left(S_{i+1/2}\Delta t\right), & u_{i+1/2} < 0.
\end{dcases}
\end{equation}
Finally, the fractional part of the Courant number is simply
\begin{equation} \label{eqn:frac_courant}
c^{frac}_{i+1/2} =
\begin{dcases}
\widehat{u}^{frac}_{i+1/2} S_{i+1/2}\Delta t / V_{i-c^{int}_{i+1/2}}, & u_{i+1/2} \geq 0, \\
\widehat{u}^{frac}_{i+1/2} S_{i+1/2}\Delta t / V_{i+1-c^{int}_{i+1/2}}, & u_{i+1/2} < 0,
\end{dcases}
\end{equation}
where the volume factor in the denominator is the volume of the cell containing the departure point.
Note that in some of the following sections, the volume factors appearing in \eqref{eqn:int_courant}, \eqref{eqn:frac_wind} and \eqref{eqn:frac_courant} will be substituted for other fields, such as a mass field.

\subsubsection{Flux Integral} \label{sec:ffsl_flux_integral}

Given the departure point $x^{dep}_{i+1/2}$, the flux of $\rho$ is computed as the integral of the mass between the flux point and the departure point as
\begin{equation}
F_{i+1/2} = 
\frac{1}{S_{i+1/2}\Delta{t}} \int^{V(x=x_{i+1/2})}_{V(x=x^{dep}_{i+1/2})} \rho \dx{V}.  
\end{equation}
As with the Courant number, the flux can be split into an integer and fractional flux, denoted $F^{int}$ and $F^{frac}$ respectively, where $F = F^{int} + F^{frac}$.
The integer part is computed by summing the integral of the field over the whole cells between the flux point and the departure point:
\begin{equation}
    F^{int}_{i+1/2} = 
    \begin{dcases}
        \frac{1}{S_{i+1/2}\Delta{t}}\sum_{j=1}^{c_{i+1/2}^{int}} \rho_{i+1-j} V_{i+1-j}, & u_{i+1/2} \geq 0, \\
        \frac{1}{S_{i+1/2}\Delta{t}}\sum_{j=1}^{|c_{i+1/2}^{int}|} \rho_{i+j} V_{i+j}, & u_{i+1/2} < 0.
    \end{dcases}
\end{equation}
The fractional part uses the subgrid reconstruction from the chosen scheme.
Although any reconstruction can be used, in this work we focus on the quadratic PPM reconstruction of $\mathcal{Q}(\xi)$ as in \eqref{eqn:quad_coef}. The fractional flux integrates the subgrid reconstruction from the facet to the departure point:
\begin{equation} \label{eqn:frac_flux}
    F^{frac}_{i+1/2} = 
    \begin{dcases}
        \frac{V_{i-c^{int}_{i+1/2}}}{S_{i+1/2}\Delta t} \int^{1}_{1-c_{i+1/2}^{frac}} \mathcal{Q}(\xi) \dx{\xi}, & u_{i+1/2} \ge 0, \\
        -\frac{V_{i+1+c^{int}_{i+1/2}}}{S_{i+1/2}\Delta t} \int^{\left| c_{i+1/2}^{frac}\right|}_{0} \mathcal{Q}(\xi) \dx{\xi}, & u_{i+1/2} < 0.  
    \end{dcases}
\end{equation}
The fractional flux of \eqref{eqn:frac_flux} can be expressed in another form, which will turn out to be key to ensuring consistent transport of tracers. 
The subgrid reconstruction $\mathcal{Q}(\xi)$ uses coefficients that are determined from field values in neighbouring grid cells.
It is therefore possible to rewrite the integral of $\mathcal{Q}(\xi)$ as a linear sum of the neighbouring field values, where the weights are polynomials in the fractional Courant $c^{frac}$, such that the fractional flux can be written as
\begin{equation} \label{eqn:consistent_flux}
F_{i+1/2}^{frac} = \widehat{u}_{i+1/2}^{frac} \mathcal{R}(\rho,c_{i+1/2}^{frac}),
\end{equation}
where $\mathcal{R}(\rho,c^{frac})$ can be interpreted as a reconstruction of $\rho$ at the facet.
Importantly, if $\rho$ is some constant $K$ then $\mathcal{R}(\rho,c_{i+1/2}^{frac})=K$ irrespective of $c_{i+1/2}^{frac}$.
Appendix A discusses the computation of $\mathcal{R}$ for PPM. \\
\\
If the density $\rho$ is a constant $K$ across the domain, then the mass flux is simply $F_{i+1/2}=Ku_{i+1/2}$, which can be seen by considering the evaluation of the integer and fractional parts of the flux.
The integer part of the flux simplifies to
\begin{equation}
    F^{int}_{i+1/2} = 
    \begin{dcases}
        \frac{K}{S_{i+1/2}\Delta{t}}\sum_{j=1}^{c_{i+1/2}^{int}}  V_{i+1-j}, & u_{i+1/2} \geq 0, \\
        \frac{K}{S_{i+1/2}\Delta{t}}\sum_{j=1}^{|c_{i+1/2}^{int}|} V_{i+j}, & u_{i+1/2} < 0,
    \end{dcases}
\end{equation}
which can be compared with \eqref{eqn:frac_flux} to show that
\begin{equation}
F^{int}_{i+1/2} = K\left(u_{i+1/2} - \widehat{u}^{frac}_{i+1/2}\right).
\end{equation}
Using \eqref{eqn:consistent_flux} with a constant $\rho$, the fractional flux is
\begin{equation}
F^{frac}_{i+1/2}=K\widehat{u}^{frac}_{i+1/2},
\end{equation}
and therefore
\begin{equation}
F_{i+1/2}=F^{int}_{i+1/2}+F_{i+1/2}^{frac} = K\left(u_{i+1/2} - \widehat{u}^{frac}_{i+1/2}\right) + K\widehat{u}^{frac}_{i+1/2} = Ku_{i+1/2}.
\end{equation}

\subsubsection{One-Dimensional Scheme for Density} 
\label{sec:1d_density}
The procedure described in Sections \ref{sec:dep_pt} and \ref{sec:ffsl_flux_integral} is used to compute departure displacements and then mass fluxes at each cell facet.
The density is then evolved in time through a simple forward-in-time scheme, so that at a given point:
\begin{equation}
\rho_i^{n+1} = \rho^n_i - \frac{\Delta t}{V_i} \left(F^x_{i+1/2}S_{i+1/2} - F^x_{i-1/2}S_{i-1/2}\right),
\end{equation}
where the superscript $n$ denotes the density at the $n$-th time level. \\
\\
For the remainder of this work, we use the operator
\begin{equation}
\mathcal{F}^x(\rho,u^{x},V)
\end{equation}
to describe the process of computing departure displacements and evaluating mass fluxes (in the $x$-direction) for the whole domain from the density $\rho$, velocity component $u^x$ and cell volume $V$, following the procedure in Sections \ref{sec:dep_pt} and \ref{sec:ffsl_flux_integral}.
The final argument represents each time that the cell volume appears in the computation of departure displacements and mass fluxes. 
With this notation, the one-dimensional scheme for density transport (in the $x$-direction) can be summarised as:
\begin{equation}
\rho^{n+1} = \rho^n - \Delta t \nabla_x\cdot \mathcal{F}^x(\rho^n,u^{x},V),
\end{equation}
where $\nabla_x\cdot$ indicates the one-dimensional discrete divergence calculation.
When $\rho^n$ is a constant $K$, the mass flux reduces to $Ku^x$ and so the value of the density at the next time level is $K(1-\Delta t\nabla_x\cdot u^x)$. This satisfies the property of divergence-linearity.

\subsubsection{One-Dimensional Scheme for an Advective Tracer} \label{sec:1d_advect}

As will be seen in the following sections, the COSMIC splitting uses a one-dimensional advective operator, which corresponds to the discretisation of \eqref{eqn:advective_transport}.
Although there is flexibility in the choice of this advective operator (for instance, \citet{lin1996ffsl} use a semi-Lagrangian scheme), in this work we take the approach of \citet{putman2007fvtransport}, using conservative operators and the transport of the unity field. \\
\\
Let $\sigma$ be the unity field, and consider the conservative transport of $\sigma m$.
Dividing the discrete conservative update of $\sigma m$ by the discrete conservative update of $\sigma$ gives a discrete advective update of $m$.
This means that the discrete advective update of $m$ in the $x$-direction is given by
\begin{equation} \label{eqn:one_d_advection_tracer_2}
m^{n+1} = m^{n} - \Delta{t} \mathcal{A}^x(m^n,u) = \left( m^n - \Delta{t} \nabla_x\cdot \mathcal{F}^x(m^n,u,V)\right)/\left( \sigma - \Delta{t} \nabla_x\cdot \mathcal{F}^x(\sigma,u,V)\right).
\end{equation}
Following \citet{putman2007fvtransport}, this can be rearranged to give a discrete operator $\mathcal{A}^x(m,u)$ described as the advective increment:
\begin{equation}
    \Delta{t} \mathcal{A}^x(m^n, u) = m^n - \left( m^n - \Delta{t} \nabla_x\cdot \mathcal{F}^x(m^n,u,V)\right)/\left( \sigma - \Delta{t} \nabla_x\cdot \mathcal{F}^x(\sigma,u,V)\right).
\end{equation}

\subsection{One-Dimensional Flux-Form Semi-Lagrangian Scheme for Conservative Tracer Transport} \label{sec:consistent_ffsl_1d}

The one-dimensional FFSL scheme is now extended to a tracer mixing ratio $m$ which evolves conservatively relative to a density $\rho$ according to \eqref{eqn:tracer_conservative}.
It is assumed throughout that $\rho$ is positive.
As discussed in Section \ref{sec:intro}, the transport of $m$ must be consistent, in the sense that a constant field is preserved.
In other words, it is required that the scheme solving for $\rho m$ produces the same results as the scheme solving for $\rho$ when $m=1$.
This is achieved by ensuring that if $m=K$ (where $K$ is a constant), the tracer flux is given by $K\bm{F}$, if the density flux is $\bm{F}$.
This property is attained by reusing the mass fluxes that are calculated during the transport of the density variable.
\\
\\
For Eulerian transport schemes, in order to provide consistency it is common to replace the transporting velocity $\bm{u}$ with the mass flux $\bm{F}$ that was used to transport the density $\rho$. This was the approach used by \citet{bendall2023solution} and is particularly relevant if $\rho$ and $m$ are not co-located.
However, care needs taking to apply this to the FFSL fluxes written in integral form.
In the approach presented here, there are two key ingredients to obtain consistency:
(a) the computation of departure points specifically for tracer transport (which are in general different from those used for the density transport), and (b) rewriting the integrals as a velocity multiplied by a facet reconstruction (see \eqref{eqn:consistent_flux}), which allows the mass flux to be used and consistency to be obtained. \\
\\
The departure points that are used for the transport of $m$ can be computed by following the approach used by \cite{skamarock2006limiters}.
These are now found by inverting the relationship
\begin{equation}
\int_{\varGamma_{i+1/2}} \bm{F} \bm{\cdot} \mathrm{d}\bm{S} = \frac{1}{\Delta t} \int_{V(x=x^{dep}_{i+1/2})}^{V(x=x_{i+1/2})} \rho \dx{V}.
\end{equation}
This means that for tracer transport, the integer and fractional components of the Courant number are defined differently to those used for the density transport, although they will coincide when $\rho$ is a constant.
The same is true of the fractional part of the flux.
However these quantities can all be computed following exactly the same procedure described in Sections \ref{sec:dep_pt} and \ref{sec:ffsl_flux_integral}, but with the following substitutions:
$m$ is used in place of $\rho$, $F^x$ takes the role of $u^x$ and the product $\rho V$ is used instead of just $V$.
For the transport of tracers, we can therefore write flux calculations through an operator of the form $\mathcal{F}^x(m,F^x,\rho V)$, whereas the equivalent operator for the transport of $\rho$ is $\mathcal{F}^x(\rho,u^x,V)$.
The fractional part of the tracer mass flux is therefore computed using the reconstruction of $m$ multiplied by the fractional part of the density flux $F^x$, and the integer part of the flux is the sum of $\rho m$ multiplied by the cell volume over the integer number of cells.
With these substitutions, if $m=K$ a constant, the tracer mass flux is $KF^x$. This follows directly from the argument at the end of Section \ref{sec:ffsl_flux_integral}. \\
\\
The overall scheme for a tracer mixing ratio for transport in the $x$-direction can then be summarised as
\begin{equation} \label{eqn:1D_density}
(\rho m)^{n+1} = (\rho m)^n - \Delta t \nabla_x\cdot \mathcal{F}^x(m^n,F^x,\rho^nV),
\end{equation}
where $(\rho m)^n\equiv \rho^n m^n$.
The mixing ratio at the next time level is then found from
\begin{equation}
m^{n+1} = \frac{(\rho m)^n}{\rho^n}.
\end{equation}
As \eqref{eqn:1D_density} is in flux form, the mass corresponding to $\rho m$ is conserved by this scheme.
If $m$ is a constant $K$, then \eqref{eqn:1D_density} becomes
\begin{equation}
(\rho m)^{n+1} = K\rho^n - \Delta t K \nabla_x\cdot F^x.
\end{equation}
As the corresponding density transport is given by
\begin{equation}
\rho^{n+1} = \rho^n - \Delta t \nabla_x \cdot F^x,
\end{equation}
the updated mixing ratio is simply $m^{n+1}=K$ and a constant mixing ratio is preserved. This satisfies the property for consistent tracer transport.

\subsection{Review of Two-dimensional COSMIC splitting} \label{sec:cosmic}

The one-dimensional FFSL methods from Section \ref{sec:1d_ffsl} can be used with splitting techniques to solve the two-dimensional transport problem. Here we review the splitting method of \citet{leonard1996cosmic} and \citet{lin1996ffsl}, known as either COSMIC or Lin-Rood splitting (here we will refer to it as COSMIC splitting). This approach has been widely used, as it provides many of the benefits of the one-dimensional methods, namely conservative large time-step-permitting transport, but in two dimensions. However, the monotonic or positive-definite limiters applied for the one-dimensional schemes in Section \ref{sec:ffsl} do not guarantee monotonicity or positivity in two dimensions (especially for large Courant number) \citep{lin1996ffsl,leonard1996cosmic}.
This will be demonstrated by numerical testing in Section \ref{sec:results2d}. \\ 
\\
The notation used will follow that from Sections \ref{sec:1d_density} and \ref{sec:1d_advect} . 
The $x$ and $y$ variables refer to local directions on a regular quadrilateral mesh, rather than the Cartesian $x$- and $y$-directions. Thus
$\mathcal{F}^x(\rho,u^x,V)$ represents the one-dimensional flux operator in the $x$-direction. 
As the fluxes are stored on cell edges, the divergence is at cell centres (co-located with $\rho$).
The advective increment operators for $\rho$ are discretisations of $\bm{u\cdot\nabla}\rho$, and are denoted by $\mathcal{A}^x(\rho,u^x)$ and $\mathcal{A}^y(\rho,u^y)$ in $x$ and $y$ respectively  (see Section \ref{sec:advective_increment} for the calculation of these increments).
Finally, the velocity fields $u^x$ and $u^y$ are given by \eqref{eqn:u_normal} (and so are defined by the mesh, rather than necessarily being physical $x$- and $y$-components), with the superscript denoting the orientation of the facet.

\subsubsection{COSMIC Splitting} \label{sec:cosmic_density}
The COSMIC splitting is applied to a density field $\rho$ as follows.
Following \citet{leonard1996cosmic} and \citet{lin1996ffsl}, an advective inner half step is performed
\begin{subequations}
\begin{align}
    \rho_A^{x/2} &= \rho^n - \frac{\Delta{t}}{2} \mathcal{A}^x(\rho^n, u^x), \\
    \rho_A^{y/2} &= \rho^n - \frac{\Delta{t}}{2} \mathcal{A}^y(\rho^n, u^y).
\end{align}
\end{subequations}
The updated field at the next time level is found from the outer conservative step, applied to the opposite inner update
\begin{equation} \label{eqn:cosmic_density}
    \rho^{n+1} = \rho^{n} - \Delta{t} \nabla_x\cdot \mathcal{F}^x(\rho_A^{y/2},u^x,V) - \Delta{t} \nabla_y \cdot \mathcal{F}^y(\rho_A^{x/2},u^y,V).
\end{equation}
The inner half steps are in advective form which provides the property of divergence-linearity, as if $\rho^n=K$, then $\rho_A^{x/2}=K$ and $\rho_A^{y/2}=K$, which gives $\rho^{n+1}=K(1-\Delta t\bm{\nabla\cdot u})$.
Hence a constant $\rho$ subject to a non-divergent flow remains constant. Using $\rho^x$ in $\mathcal{F}^y$ and $\rho^y$ in $\mathcal{F}^x$ ensures that all of the cross-terms are captured \citep{kent2019crossterms} and that the scheme is symmetric. 

\subsubsection{Advective Increment} \label{sec:advective_increment}

The one-dimensional advective form increment of the density field can be computed using the schemes in Section \ref{sec:ffsl} following the method given in Section \ref{sec:1d_advect}. The advective increments are thus
\begin{subequations} \label{eqn:advective_operators}
\begin{align}
    \Delta{t} \mathcal{A}^x(\rho^n, u^x) &= \rho^n - \left( \rho^n - \Delta{t} \nabla_x\cdot \mathcal{F}^x(\rho^n,u^x,V)\right)/\left( \sigma - \Delta{t} \nabla_x\cdot \mathcal{F}^x(\sigma,u^x,V)\right), \\
    \Delta{t} \mathcal{A}^y(\rho^n, u^y) &= \rho^n - \left( \rho^n - \Delta{t} \nabla_y\cdot \mathcal{F}^y(\rho^n,u^y,V)\right)/\left( \sigma - \Delta{t} \nabla_y\cdot \mathcal{F}^y(\sigma,u^y,V)\right),
\end{align}
\end{subequations}
where $\sigma$ is the unity field.

\subsubsection{Consistent COSMIC Splitting for Conservative Tracers} \label{sec:cosmic_tracer}
The COSMIC scheme can be made consistent for tracer transport following \cite{lin1996ffsl}.
This concerns a tracer mixing ratio $m$ and the dry density $\rho$, with corresponding dry mass fluxes $F^x$ and $F^y$, which are computed from \eqref{eqn:cosmic_density}. \\ 
\\
The inner step is performed on the mixing ratio only
\begin{subequations}
\begin{align}
    m^{x/2} &= m^n - \frac{\Delta{t}}{2} \mathcal{A}^x(m^n, u^x), \\
    m^{y/2} &= m^n - \frac{\Delta{t}}{2} \mathcal{A}^y(m^n, u^y).
\end{align}
\end{subequations}

The outer step now uses the dry mass fluxes in place of the velocities to compute the tracer fluxes.
Following \citet{skamarock2006limiters}, new departure points and fractional mass fluxes need computing, according to the procedure outlined in Section \ref{sec:consistent_ffsl_1d}.
Then the updated mixing ratio is found from
\begin{equation} \label{eqn:cosmic_tracer_end}
   (\rho m)^{n+1} = (\rho m)^{n} - \Delta{t} \nabla_x\cdot \mathcal{F}^x(m^{y/2},F^x,\rho^nV) - \Delta{t} \nabla_y\cdot \mathcal{F}^y(m^{x/2},F^y,\rho^nV),
\end{equation}
where $m^{n+1}$ can be backed out by dividing by $\rho^{n+1}$. This approach to consistent COSMIC splitting is used for the numerical results presented in this paper.

\section{Two-dimensional Monotonic SWIFT Splitting} \label{sec:swift_splitting} 

The COSMIC splitting from Section \ref{sec:cosmic} provides large time-step-permitting conservative and consistent transport. However, it does not guarantee monotonicity when monotonic limiters are used for the one-dimensional schemes, especially for large Courant number flow (see \citep{leonard1996cosmic,lin1996ffsl,bott2010improving} and the results from Section \ref{sec:results2d}). \\
\\
Here we present a new splitting formulation as part of the SWIFT scheme. This splitting can be used for conservative or advective transport, ensures consistency of tracers, is large time-step-permitting, and enforces monotonicity or positivity in two-dimensions when monotonic or positive-definite limiters are used for the one-dimensional scheme. This means that the monotonic or positivity-preserving properties of the one-dimensional schemes are maintained for the two-dimensional method. 

\subsection{SWIFT Splitting for Density} \label{sec:swift_density}
This section describes the application of the SWIFT splitting to density transport.
Let $\sigma$ be the unity field and $\rho$ the dry density field.
The inner step for density transport consists of conservative steps:
\begin{subequations} \label{eqn:swift_density_inner}
\begin{align} \label{eqn:swift_density_x}
    \rho_I^x &= \rho^n - \Delta{t} \nabla_x\cdot \mathcal{F}^x(\rho^n, u^x,V), \\
    \rho_I^y &= \rho^n - \Delta{t} \nabla_y\cdot \mathcal{F}^y(\rho^n, u^y,V).
\end{align}
\end{subequations}
The unity field is also transported in conservative form, so that
\begin{subequations} \label{eqn:sigma_x_sigma_y}
\begin{align}
    \sigma^x &= \sigma - \Delta{t} \nabla_x\cdot \mathcal{F}^x(\sigma, u^x,V), \\
    \sigma^y &= \sigma - \Delta{t} \nabla_y\cdot \mathcal{F}^y(\sigma, u^y,V).
\end{align}
\end{subequations}
Advective updates are defined using the transported unity fields:
\begin{subequations}
\begin{align} \label{eqn:swift_density_adv_x}
    \rho^x_A &= \rho_I^x / \sigma^x, \\ \label{eqn:swift_density_adv_y}
    \rho^y_A &= \rho_I^y / \sigma^y.
\end{align}
\end{subequations}
Densities with subscript $I$ denote values obtained from the inner step, while this with subscript $A$ include the division by an appropriately transported unity field to give a density with an advective update.
The departure points and fractional components of the wind are recomputed using $\sigma^yV$ and $\sigma^xV$ as the \textit{volume} in equations \eqref{eqn:int_courant}-\eqref{eqn:frac_courant}.
The outer density step computes the fluxes using reconstructions of $\rho^x_A$ and $\rho^y_A$ as in the COSMIC splitting, but now updates the fields from the inner steps:
\begin{equation} \label{eqn:swift_density}
\rho^{n+1} = \frac{1}{2}\left[\rho_I^{y} - \Delta{t}\nabla_x\cdot \mathcal{F}^x(\rho^y_A, u^x, \sigma^yV) + \rho_I^x - \Delta{t}\nabla_y\cdot \mathcal{F}^y(\rho^x_A, u^y, \sigma^xV)\right].
\end{equation}
This essentially combines the fluxes from the inner and outer steps, as
\begin{equation} \label{eqn:swift_density_2}
\begin{split}
\rho^{n+1} \equiv \rho^n & - \frac{\Delta t}{2} \nabla_x\cdot \mathcal{F}^x(\rho^n, u^x, V) - \frac{\Delta{t}}{2} \nabla_x\cdot \mathcal{F}^x(\rho^y_A, u^x, \sigma^yV) \\
& - \frac{\Delta{t}}{2} \nabla_y\cdot \mathcal{F}^y(\rho^n, u^y,V) - \frac{\Delta{t}}{2} \nabla_y\cdot \mathcal{F}^y(\rho^x_A, u^y, \sigma^xV).
\end{split}
\end{equation}
The one-dimensional density updates are stored, as are the dry mass fluxes, as both are needed for the SWIFT tracer transport:
\begin{subequations}
\begin{align}
F^x &= \frac{1}{2}\left[ \mathcal{F}^x(\rho^n, u^x,V) +  \mathcal{F}^x(\rho^y_A, u^x, \sigma^yV)\right], \\ 
F^y &= \frac{1}{2}\left[ \mathcal{F}^y(\rho^n, u^y,V) +  \mathcal{F}^y(\rho^x_A, u^y, \sigma^xV)\right], \\
\rho^x &= \rho^{n} - \Delta{t} \nabla_x\cdot F^x, \\
\rho^y &= \rho^{n} - \Delta{t} \nabla_y\cdot F^y. 
\end{align}
\end{subequations}

\subsection{SWIFT Splitting for Consistent Tracers} \label{sec:swift_tracer}
The SWIFT splitting is now applied to tracer transport.
This splitting is designed to preserve the consistency, monotonicity and positivity properties of the one-dimensional calculation, whilst still providing local conservation of the tracer and stability that is not conditional on the Courant number.
The steps for SWIFT splitting for tracers are as follows. \\
\\
First, the departure points and dry mass fluxes are recomputed using $\rho^nV$ and the dry mass fluxes substituted into equations \eqref{eqn:int_courant}-\eqref{eqn:frac_courant} for the first transport step.
The inner transport step for the tracer is then given by:
\begin{subequations}
\begin{align}  \label{eqn:swift_split_1}   
        (\rho m)^x &= (\rho m)^n - \Delta{t} \nabla_x\cdot \mathcal{F}^x(m^n,F^{x},\rho^nV), \\ \label{eqn:swift_split_2}
        (\rho m)^y &= (\rho m)^n - \Delta{t} \nabla_y\cdot \mathcal{F}^y(m^n,F^{y},\rho^nV).
\end{align}
\end{subequations}
Updated mixing ratio fields are found by dividing by the conservative densities following one-dimensional transport steps
\begin{subequations}
\begin{align}
m^x &= (\rho m)^x / \rho^x, \\
m^y &= (\rho m)^y / \rho^y.
\end{align}
\end{subequations}
The departure points and fractional mass fluxes are recomputed using $\rho^y$ and $\rho^x$ for the second transport step.
Then the updated tracer density combines fluxes from the inner and outer steps:
\begin{equation} \label{eqn:swift_tracer}
\begin{split}
(\rho m)^{n+1} = (\rho m)^n & -\frac{\Delta{t}}{2}\nabla_x\cdot \mathcal{F}^x(m^n,F^{x},\rho^nV) -\frac{\Delta{t}}{2}\nabla_x\cdot \mathcal{F}^x(m^y,F^{x},\rho^yV) \\
        & - \frac{\Delta{t}}{2}\nabla_y\cdot \mathcal{F}^y(m^n,F^{y},\rho^nV) -\frac{\Delta{t}}{2}\nabla_y\cdot \mathcal{F}^y(m^x,F^{y},\rho^xV).
\end{split}
\end{equation}
This can also be written as
\begin{equation}
(\rho m)^{n+1} = \frac{1}{2}\left[(\rho m)^y - \Delta t \nabla_x\cdot \mathcal{F}^x(m^y,F^{x},\rho^y, V) + (\rho m)^x - \Delta t \nabla_y\cdot \mathcal{F}^y(m^x,F^{y},\rho^x, V) \right].
\end{equation}
The tracer mixing ratio at the end of the time step is then computed from the tracer density and the dry density
\begin{equation}
    m^{n+1} = (\rho m)^{n+1} / \rho^{n+1}.
\end{equation}

\subsection{SWIFT Splitting for Advective Tracers}

The SWIFT splitting can also be applied to a mixing ratio that is solved advectively, in the form of \eqref{eqn:advective_transport}. Note that mixing ratios solved this way will not be conservative or consistent with the dry density. \\ 
\\
For the advective form of SWIFT splitting we follow the steps given in Section \ref{sec:swift_density} but replace $\rho$ with the mixing ratio $m$. The final step is to then divide through by unity transported for a whole time step
\begin{equation}
\begin{split}
(\sigma m)^{n+1} = m^{n} & - \frac{\Delta{t}}{2}\nabla_x\cdot \mathcal{F}^x(m^n, u^x, V) - \frac{\Delta{t}}{2}\nabla_x\cdot \mathcal{F}^x(m^y_A, u^x, \sigma^yV) \\
& - \frac{\Delta{t}}{2}\nabla_y\cdot \mathcal{F}^y(m^n, u^y, V) - \frac{\Delta{t}}{2}\nabla_y\cdot \mathcal{F}^y(m^x_A, u^y, \sigma^xV),
\end{split}
\end{equation}
where
\begin{align}
    m^x_A &= (m^n - \Delta{t} \nabla_x\cdot \mathcal{F}^x(m^n, u^x, V)) / \sigma^x, \\
    m^y_A &= (m^n - \Delta{t} \nabla_y\cdot \mathcal{F}^y(m^n, u^y, V)) / \sigma^y.
\end{align}
The updated mixing ratio is then found from
\begin{equation}
m^{n+1} = (\sigma m)^{n+1}/\sigma^{n+1},
\end{equation}
where $\sigma^{n+1}$ is calculated from
\begin{equation}
\sigma^{n+1} = \sigma - \Delta{t} \nabla_x\cdot \mathcal{F}^x(\sigma, u^x, V) -  \nabla_y \cdot \mathcal{F}^y(\sigma, u^y, V).
\end{equation}
This approach is also equivalent to the conservative tracer transport scheme of Section \ref{sec:swift_tracer}, but with $\rho$ replaced by $\sigma$.

\subsection{Analysis of SWIFT Splitting}

This section demonstrates how the SWIFT splitting ensures consistency between tracers and density, how the monotonicity or positivity properties of the one-dimensional scheme are achieved for the two-dimensional splitting, and how the SWIFT splitting is similar to other well-known splitting techniques. 

\subsubsection{Divergence-Linearity}
The SWIFT splitting for density transport presented in Section \ref{sec:swift_density}
provides increments that are linear in the divergence.
This can be seen by inserting $\rho^n=K$ for constant $K$ into \eqref{eqn:swift_density_inner}, which when combined with \eqref{eqn:sigma_x_sigma_y} returns densities for the outer step of $\rho_A^x=K$ and $\rho_A^y=K$.
These values can then be substituted into \eqref{eqn:swift_density}.
Using the property discussed in Section \ref{sec:ffsl_flux_integral} that when the density is constant, the mass flux is proportional to the wind field, the resulting density is
\begin{equation}
\rho^{n+1} = K(1-\Delta t \bm{\nabla\cdot u}).
\end{equation}
Thus, the increment to the field is linear in the divergence, and a constant density is preserved under divergence-free flow.

\subsubsection{Consistency}

The SWIFT scheme is consistent as it preserves a constant mixing ratio. Let $m=K$ be constant and assume that the dry density $\rho$ is positive-definite. The inner stage horizontal fluxes are equal to the constant multiplied by the dry flux
\begin{subequations}
\begin{align}   
   \mathcal{F}^x(m^n,F^{x},\rho^nV) &= K F^x, \\
   \mathcal{F}^y(m^n,F^{y},\rho^nV) &= K F^y,
\end{align}
\end{subequations}
so the first stage updates are
\begin{subequations}
\begin{align}        
    (\rho m)^x &=\rho^n K - \Delta{t} \nabla_x\cdot \mathcal{F}^x(K,F^{x},\rho^nV) = \rho^n K - \Delta{t} K \nabla_x\cdot F^x = K \rho^x, \\
    (\rho m)^y &=\rho^n K - \Delta{t} \nabla_y\cdot \mathcal{F}^y(K,F^{y},\rho^nV) = \rho^n K - \Delta{t} K \nabla_y\cdot F^y = K \rho^y, \\
    \Rightarrow m^x &= (\rho m)^x / \rho^x = K, \\
    \Rightarrow m^y &= (\rho m)^y / \rho^y = K.
\end{align}
\end{subequations}
The outer stage fluxes, using the recomputed dry mass fluxes from $\rho^y$ and $\rho^x$, are
\begin{subequations}
\begin{align}
    \mathcal{F}^x(m^y,F^{x},\rho^yV) & = K F^x, \\
    \mathcal{F}^y(m^x,F^{y},\rho^xV) & = K F^y,
\end{align}
\end{subequations}
 and hence the update of the tracer is
 \begin{subequations}
\begin{align}
    (\rho m)^{n+1} &= \rho^n K -\frac{\Delta{t}}{2}\nabla_x\cdot (K F^x) -\frac{\Delta{t}}{2}\nabla_x\cdot (K F^x) -\frac{\Delta{t}}{2}\nabla_y\cdot (K F^y) -\frac{\Delta{t}}{2}\nabla_y\cdot (K F^y) = K \rho^{n+1}, \\
    \Rightarrow m^{n+1} &=  (\rho m)^{n+1}/ \rho^{n+1} = K,
\end{align}
\end{subequations}
and the constant tracer is preserved.

\subsubsection{Multi-dimensional Monotonicity Properties of One-Dimensional Schemes}
Consider the one-dimensional schemes in Section \ref{sec:ffsl} with the application of a monotonic or positive-definite limiter. The one-dimensional solution is monotonic or positive-definite provided the input to the fluxes is the same as the field being transported. For a consistent tracer transport example, a monotonic scheme used to compute the fluxes in
\begin{equation} \label{eqn:mono_condition}
    (\rho m)^{n+1} = (\rho m)^n - \Delta{t} \nabla_x\cdot \mathcal{F}^x(m^n, F^x, \rho^nV), 
\end{equation}
will produce a shape-preserving solution for $(\rho m)$, and hence a monotonic solution for $m$. The same holds for a positive-definite limiter and positivity.
Crucially, the arguments $m^n$ and $\rho^n$ to the flux operator $\mathcal{F}^x$ must be the same mixing ratio and density that are being updated. \\
\\
The COSMIC splitting is not monotonic for large Courant numbers due to the fluxes in \eqref{eqn:cosmic_tracer_end} acting on $m^y$ and $m^x$, while the field that is transported is $(\rho m)^n$. This is not in the same form as \eqref{eqn:mono_condition}. \\ 
\\
The SWIFT splitting is designed to ensure multi-dimensional monotonicity when one-dimensional limiters are used. Consider the inner SWIFT step. Equations \eqref{eqn:swift_split_1} and \eqref{eqn:swift_split_2} are both in the form of \eqref{eqn:mono_condition}, hence $(\rho m)^x$ and $(\rho m)^y$, and thus $m^x$ and $m^y$ are monotonic. \\
\\
From \eqref{eqn:swift_tracer} the outer SWIFT step is written in the required form:
\begin{equation}
    (\rho m)^{n+1}  = \frac{1}{2}(\rho m)^x  -\frac{\Delta{t}}{2}\nabla_y\cdot \mathcal{F}^y(m^x,F^{y},\rho^xV) +  \frac{1}{2}(\rho m)^y - \frac{\Delta{t}}{2}\nabla_x\cdot \mathcal{F}^x(m^y,F^{x},\rho^yV). \label{eqn:swift_is_monotonic}
\end{equation}
As the departure points in the outer step have been recomputed using $\rho^x$ and $\rho^y$, this is in the form of \eqref{eqn:mono_condition}.
Equation \eqref{eqn:swift_is_monotonic} shows that the final mixing ratio field comes from the average of two updates of the form \eqref{eqn:mono_condition}.
The averaging does not generate new extrema, so any monotonicity or positive-definite properties are preserved in the two-dimensional splitting.

\subsubsection{Stability and Relationship to COSMIC}
The inner step fields $\rho_A^{x/2}$ and $\rho_A^{x}$ of the COSMIC and SWIFT splittings can be related to one another through
\begin{equation}
\rho_A^{x/2} = \frac{1}{2}\rho^n + \frac{1}{2}\rho_A^x,
\end{equation}
and similarly for the step in the $y$-direction.
If no monotonic or positive-definite limiter is used, the flux operator $\mathcal{F}$ is linear in its first argument and the COSMIC splitting from \eqref{eqn:cosmic_density} can be rewritten in terms of the quantities used in the SWIFT splitting to give:
\begin{equation} \label{eqn:cosmic_linear}
\begin{split}
\rho^{n+1} = \rho^n & - \frac{\Delta t}{2} \nabla_x\cdot \mathcal{F}^x(\rho^n, u^x,V) - \frac{\Delta{t}}{2} \nabla_x\cdot \mathcal{F}^x(\rho^y_A, u^x,V) \\
& - \frac{\Delta{t}}{2} \nabla_y\cdot \mathcal{F}^y(\rho^n, u^y,V) - \frac{\Delta{t}}{2} \nabla_y\cdot \mathcal{F}^y(\rho^x_A, u^y,V).
\end{split}
\end{equation}
This can be compared directly with the transport of density using SWIFT, and so repeating \eqref{eqn:swift_density_2} here
\begin{equation}
\begin{split}
\rho^{n+1} \equiv \rho^n & - \frac{\Delta t}{2} \nabla_x\cdot \mathcal{F}^x(\rho^n, u^x, V) - \frac{\Delta{t}}{2} \nabla_x\cdot \mathcal{F}^x(\rho^y_A, u^x, \sigma^yV) \\
& - \frac{\Delta{t}}{2} \nabla_y\cdot \mathcal{F}^y(\rho^n, u^y,V) - \frac{\Delta{t}}{2} \nabla_y\cdot \mathcal{F}^y(\rho^x_A, u^y, \sigma^xV),
\end{split}
\end{equation}
reveals that, in the absence of a limiter, the COSMIC and SWIFT splittings differ through the inclusion of $\sigma^x$ and $\sigma^y$ arguments for the outer flux evaluations and thus the computation of the departure points for the outer step.
\\
\\
A corollary of this is that if the transporting velocity field is constant in space, then the COSMIC splitting of Section \ref{sec:cosmic_density} and the SWIFT splitting of Section \ref{sec:swift_density} for density are equivalent.
This can be seen by noting that $\sigma^x=\sigma^y=1$ if the flow is constant. \\
\\
A consequence of the equivalence of these schemes under constant flow is that the stability analysis of SWIFT for density fields also reduces to that of COSMIC, since standard stability analysis requires the scheme to be applied to a linear problem.
A thorough stability analysis of the COSMIC splitting using PPM was performed by \citet{lauritzen2007fvstab}, so this same analysis applies to the SWIFT splitting.
As demonstrated by \citet{lin1996ffsl}, the splitting is stable provided the Lipschitz condition is not broken in either of the horizontal directions. \\
\\
For tracer transport, linear stability analysis would also require the underlying density field to be constant.
With this additional assumption of constant density, along with the constant velocity, the SWIFT splitting for tracer transport also reduces to the COSMIC splitting for tracer transport (provided the one-dimensional FFSL scheme is linear).
With constant velocity, an initially constant $\rho$ will remain constant, and then both the COSMIC and SWIFT tracer splittings are equivalent, and also reduce to the splittings used for density transport. \\
\\
Another approach to dimensional splitting is a Strang splitting of the $x$ and $y$ calculations as described by \citet{skamarock2006limiters}, which also creates a two-dimensional scheme which retains any monotonicity or positivity properties of the one-dimensional scheme.
However, such a splitting breaks the symmetry between the two horizontal dimensions that is inherent in the underlying model.
Although a Strang splitting involves fewer flux computations than the SWIFT splitting, it offers no advantage in terms of parallel data communications.
In an MPI parallelisation with a domain partitioning, both splittings require three halo exchanges.

\section{Two-Dimensional Results} \label{sec:results2d}

We test the SWIFT splitting by running a number of prescribed velocity transport tests on the doubly-periodic plane with $x$ and $y$ ranging from $-0.5$ km to $0.5$ km, and comparing with COSMIC splitting. For each test we use PPM as the one-dimensional scheme, and where noted use the strict limiting. For each test the density $\rho$ is computed using the same splitting as the tracer $m$.
All tracers are transported conservatively, have an initial minimum of 0 and maximum of 1.
The setup of the test cases is given in Appendix B, but in all tests the scalar fields are transported once around the domain so that the analytic final solutions are equal to their initial values.
Throughout the results it is demonstrated that at small Courant numbers, the behaviour of the SWIFT and COSMIC splittings are very similar and so SWIFT retains the benefits of the COSMIC splitting. However at large Courant numbers, the COSMIC splitting does not maintain monotonicity when a limiter is used, while the SWIFT splitting does.

\subsection{Constant Velocity Test} \label{sec:test1}

\begin{figure}[h!]
\centering
\includegraphics[width=0.9\textwidth]{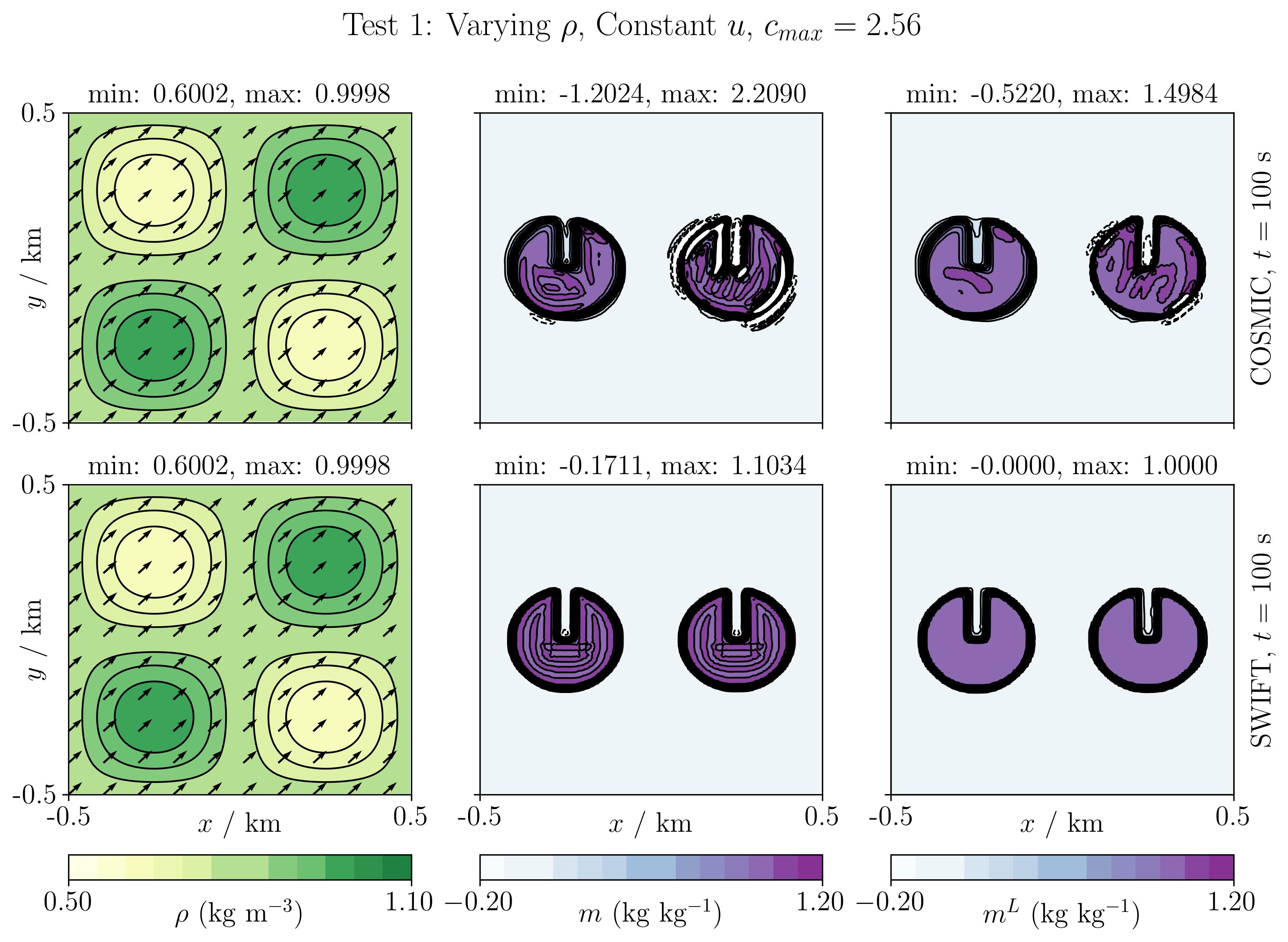}
\caption{The final transported fields from the constant velocity test of Section \ref{sec:test1}, using a spatially-varying density. These fields are taken from the large Courant number case, with $\Delta t=2$ s. The top row shows fields transported with the COSMIC splitting, while the bottom row contains fields transported with the SWIFT splitting.
For the density field $\rho$ (left column) the two splittings give the same result.
The central column displays a conservative mixing ratio, with no limiter applied.
The limited tracer (right column) is monotonic with the SWIFT splitting but not with the COSMIC splitting.
The density fields are contoured at intervals of $0.05$ kg m$^{-3}$, with the $0.8$ kg m$^{-3}$ contour omitted.
The mixing ratio fields have contours every $0.1$ kg kg$^{-1}$, omitting the zero contour.
Arrows on the density plots show the direction and magnitude of the transporting velocity.}
\label{fig:test1}
\end{figure}

The constant velocity test is used with both an initially-constant density, and a spatially-varying density.
The tracers are initialised as slotted cylinders on the $128\times 128$ grid.
The test is run with two different time steps: $0.2$ s which gives a maximum one-dimensional Courant number of $c_{max}=0.256$, and $2$ s which corresponds to $c_{max}=2.56$.
\\
\begin{table}[h!]
\small
\begin{center}
\begin{tabular}{| c | l | c c c | c c c | } \hline
  & & \multicolumn{3}{c|}{$c_{max}=0.256$} & \multicolumn{3}{c|}{$c_{max}=2.56$} \\
  & Scheme &  Min & Max & $L^2$ error & Min & Max & $L^2$ error   \\ \hline
  \parbox[t]{2mm}{\multirow{4}{*}{\rotatebox[origin=c]{90}{Const. $\rho$}}}
  & COSMIC $m$ & -0.116 & 1.233 & 2.21$\times 10^{-1}$ & -0.197 & 1.119 & 1.74$\times 10^{-1}$ \\
  & SWIFT $m$ & -0.116 & 1.233 & 2.21$\times 10^{-1}$ & -0.197 & 1.119 & 1.74$\times 10^{-1}$ \\
  & COSMIC $m^L$ & 0.000 & 0.998 & 2.55$\times 10^{-1}$ & -8.0$\times 10^{-3}$ & 1.004 & 1.81$\times 10^{-1}$ \\
  & SWIFT  $m^L$ & 0.000 & 1.000 & 2.53$\times 10^{-1}$ & 0.000 & 1.000 & 1.87$\times 10^{-1}$ \\ \hline
  \parbox[t]{2mm}{\multirow{6}{*}{\rotatebox[origin=c]{90}{Varying $\rho$}}}
  & COSMIC $\rho$ & 0.600 & 1.000 & 1.10$\times 10^{-6}$ & 0.600 & 1.000 & 1.83$\times 10^{-7}$ \\
  & SWIFT $\rho$ & 0.600 & 1.000 & 1.10$\times 10^{-6}$ & 0.600 & 1.000 & 1.83$\times 10^{-7}$ \\
  & COSMIC $m$ & -0.119 & 1.239 & 2.21$\times 10^{-1}$ & -1.261 & 2.283 & 3.15$\times 10^{-1}$ \\
  & SWIFT $m$ & -0.118 & 1.236 & 2.21$\times 10^{-1}$ & -0.191 & 1.121 & 1.76$\times 10^{-1}$ \\
  & COSMIC $m^L$ & 0.000 & 0.998 & 2.54$\times 10^{-1}$ & -0.469 & 1.438 & 2.19$\times 10^{-1}$ \\
  & SWIFT $m^L$ & 0.000 & 0.998 & 2.54$\times 10^{-1}$ & 0.000 & 1.000 & 1.88$\times 10^{-1}$ \\ \hline
\end{tabular}
\caption{Results for the constant velocity tests. The top set are for a constant density field, the bottom are for the varying density. The normalised $L^2$ errors, minimum and maximum values are given for slotted cylinder initial conditions on the $128\times128$ grid, for Courant numbers $c_{max}=0.256$ and $c_{max}=2.56$.
The mixing ratio $m$ has no limiter applied, while $m^L$ uses the limiter of Appendix A\ref{sec:limiter}.}
\label{table:test1}
\end{center}
\end{table}
\\
Statistics from the end of these simulations (at $t=100$ s) are shown in Table \ref{table:test1}.
The ``Min'' and ``Max'' columns give the instantaneous minima and maxima of the transported fields at the end of the final time step.
The tracer field has an initial minimum of 0 and an initial maximum of 1, while the spatially-varying density has an initial minimum of 0.6 and an initial maximum of 1.
The ``$L^2$ error'' column describes the normalised $L^2$ error, which is given by $||q - q_{true}||/||q_{true}||$ for transported field $q$ and analytic solution $q_{true}$.
The left-hand side of the table presents the values from a run with the smaller Courant number, while the right-hand side contains the values for the larger Courant number case.
Each row represents a different transported variable.
The superscript for $m^L$ indicates that the monotonic limiter was used for the transport of that mixing ratio.
The upper half of the table contains statistics for tracers transported relative to an initially-constant density, while the lower half describes the case with a spatially-varying density.
The final fields corresponding to the lower-right corner of the table are plotted in Figure \ref{fig:test1}.\\
\\
The results in the top half of Table \ref{table:test1} demonstrate that for constant flow and constant density, the COSMIC and SWIFT splittings are equivalent in the absence of a limiter. The results in the bottom half of Table \ref{table:test1} demonstrate that the COSMIC and SWIFT splittings are equivalent for density under constant flow, but not for tracers when the density varies. 
For the small Courant number tests, the limited tracer with COSMIC splitting does not generate extrema beyond the initial extrema.
However for large Courant number tests, the limited tracer is no longer bounded by its initial extrema.
In both cases, the limited tracer with the SWIFT splitting remains limited by its initial extrema.
For all variables, the SWIFT error norms are of similar magnitude to those with the COSMIC splitting.

\subsection{Non-Divergent Deformational Test} \label{sec:test2}
We now use the deformational velocity of \cite{skamarock2006limiters} with background flow (see \citet{kent2020positive}), and slotted cylinder initial tracers. As before, we consider both cases with small and large Courant numbers, and a constant and a spatially-varying density. It is worth noting that there are very few tracer transport test cases with a varying density in the literature, and so this test, with non-constant density, deformational flow, and large Courant number, is a very challenging test of the numerical tracer transport scheme. 
The two different time steps were again $\Delta t=0.2$ s and $2$ s, which gave maximum one-dimensional Courant numbers of $0.512$ and $5.12$ respectively. \\
\begin{figure}[ht!]
\centering
\includegraphics[width=0.9\textwidth]{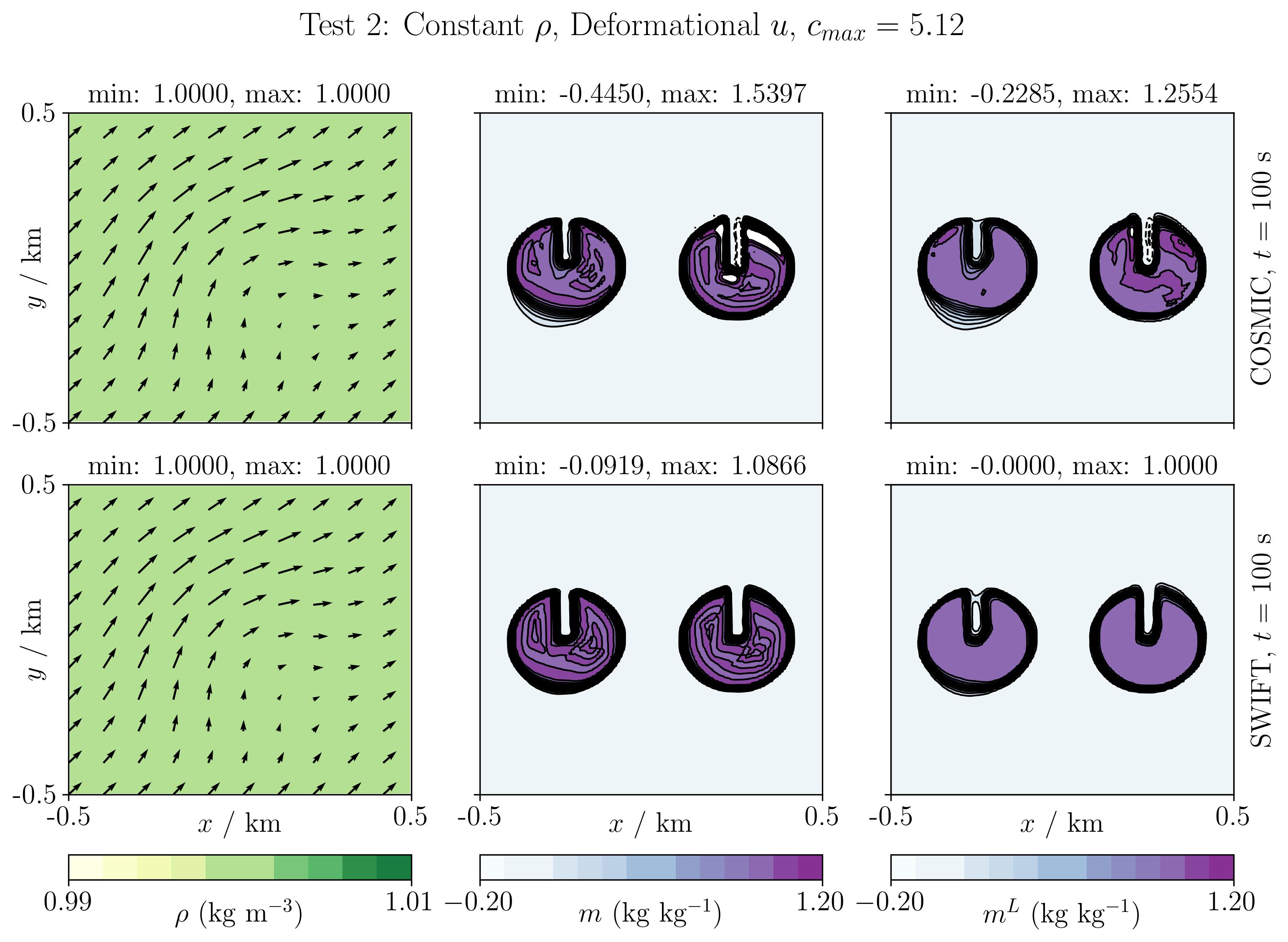}
\caption{The transported fields at the end of the final time step ($t=100$ s) of the non-divergent deformational test of Section \ref{sec:test2} with a constant density. These fields are taken from the large Courant number case, with $\Delta t=2$ s. Different rows show the fields transported with the COSMIC splitting and SWIFT splittings.
The mixing ratio fields have contours every $0.1$ kg kg$^{-1}$, omitting the zero contour.
Arrows on the density plots show the direction and magnitude of the transporting velocity.}\label{fig:test2constrho}
\end{figure} \\
\\
Figure \ref{fig:test2constrho} shows the fields at the end of the final time step for the constant density case with maximum one-dimensionl Courant number $5.12$ on the $128\times128$ grid. As the velocity is non-divergent the density remains constant at all times for both COSMIC and SWIFT. The plots for COSMIC and SWIFT are very similar for the tracer, with the shape of the cylinder being maintained by both schemes. With the limiter applied, COSMIC still produces over-and under-shoots, although they are small in magnitude, whereas SWIFT ensures the tracer is bounded by its initial values. This test demonstrates that both COSMIC and SWIFT produce accurate results when using PPM, even for large Courant numbers, but only SWIFT ensures monotonicity. \\
\begin{figure}[ht!]
\centering
\includegraphics[width=0.78\textwidth]{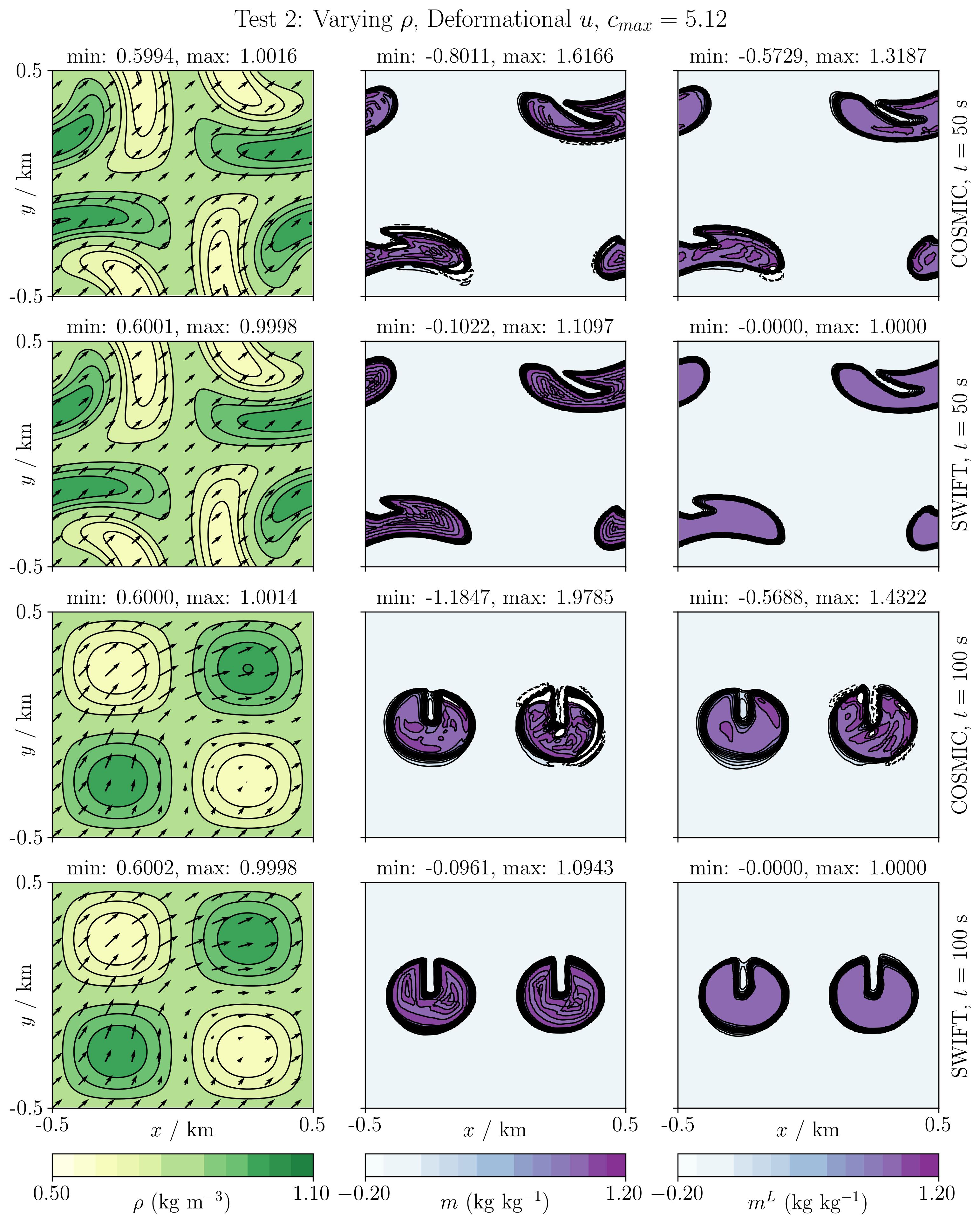}
\caption{The transported fields at the halfway point ($t=50$ s, top two rows) and at the end of the final time step ($t=100$ s, bottom two rows) of the non-divergent deformational test of Section \ref{sec:test2}. These fields are taken from the large Courant number case, with $\Delta t=2$ s. Different rows show the fields transported with the COSMIC splitting and SWIFT splittings.
Even under strong deformation, the limited tracer with the SWIFT splitting remains monotonic.
The density fields are contoured at intervals of $0.05$ kg m$^{-3}$, with the $0.8$ kg m$^{-3}$ contour omitted.
The mixing ratio fields have contours every $0.1$ kg kg$^{-1}$, omitting the zero contour.
Arrows on the density plots show the direction and magnitude of the transporting velocity.}\label{fig:test2}
\end{figure} 

Results for the spatially-varying density case are presented in Table \ref{table:test2}, which has the same structure as Table \ref{table:test1}. 
This again demonstrates that transport with the SWIFT splitting is monotonic when the COSMIC splitting is not, whilst error norms are still comparable.
Figure \ref{fig:test2} plots the fields for the larger Courant number case, at the halfway point and at the end of the simulation. Both COSMIC and SWIFT splitting produce accurate solutions for the density. For the tracer COSMIC splitting produces large over- and under-shoots, even when the limiter is used. The SWIFT splitting is monotonic if the limiter is used, but even for unlimited PPM the over- and under-shoots are smaller than for the COSMIC splitting. \\
\begin{table}[h!]
\small
\begin{center}
\begin{tabular}{| c | l | c c c | c c c |} \hline
  &  & \multicolumn{3}{c|}{$c_{max}=0.512$} & \multicolumn{3}{c|}{$c_{max}=5.12$} \\
  & Scheme &  Min & Max & $L^2$ error & Min & Max & $L^2$ error   \\ \hline
  \parbox[t]{2mm}{\multirow{6}{*}{\rotatebox[origin=c]{90}{Varying $\rho$}}}
  & COSMIC $\rho$ & 0.600 & 1.000 & 2.26$\times 10^{-5}$ & 0.600 & 1.001 & 1.68$\times 10^{-3}$ \\
  & SWIFT $\rho$ & 0.600 & 1.000 & 1.94$\times 10^{-5}$ & 0.600 & 1.000 & 1.37$\times 10^{-3}$ \\
  & COSMIC $m$ & -0.170 & 1.206 & 2.37$\times 10^{-1}$ & -1.228 & 1.978 & 3.07$\times 10^{-1}$ \\
  & SWIFT $m$ & -0.167 & 1.208 & 2.36$\times 10^{-1}$ & -0.108 & 1.102 & 1.84$\times 10^{-1}$ \\
  & COSMIC $m^L$ & 0.000 & 0.998 & 2.67$\times 10^{-1}$ & -0.510 & 1.427 & 2.56$\times 10^{-1}$ \\
  & SWIFT $m^L$ & 0.000 & 0.998 & 2.66$\times 10^{-1}$ & 0.000 & 1.000 & 2.08$\times 10^{-1}$ \\ \hline
\end{tabular}
\caption{Statistics from the non-divergent deformational test. The normalised $L^2$ errors, minimum and maximum values are given for the slotted cylinder initial conditions on the $128\times128$ grid, for Courant numbers $c_{max}=0.512$ and $c_{max}=5.12$.
The mixing ratio $m$ has no limiter applied, while $m^L$ uses the limiter of Appendix A\ref{sec:limiter}.}
\label{table:test2}
\end{center}
\end{table}

\subsection{Divergent Deformational Test} \label{sec:test3}

\begin{figure}[ht!]
\centering
\includegraphics[width=0.78\textwidth]{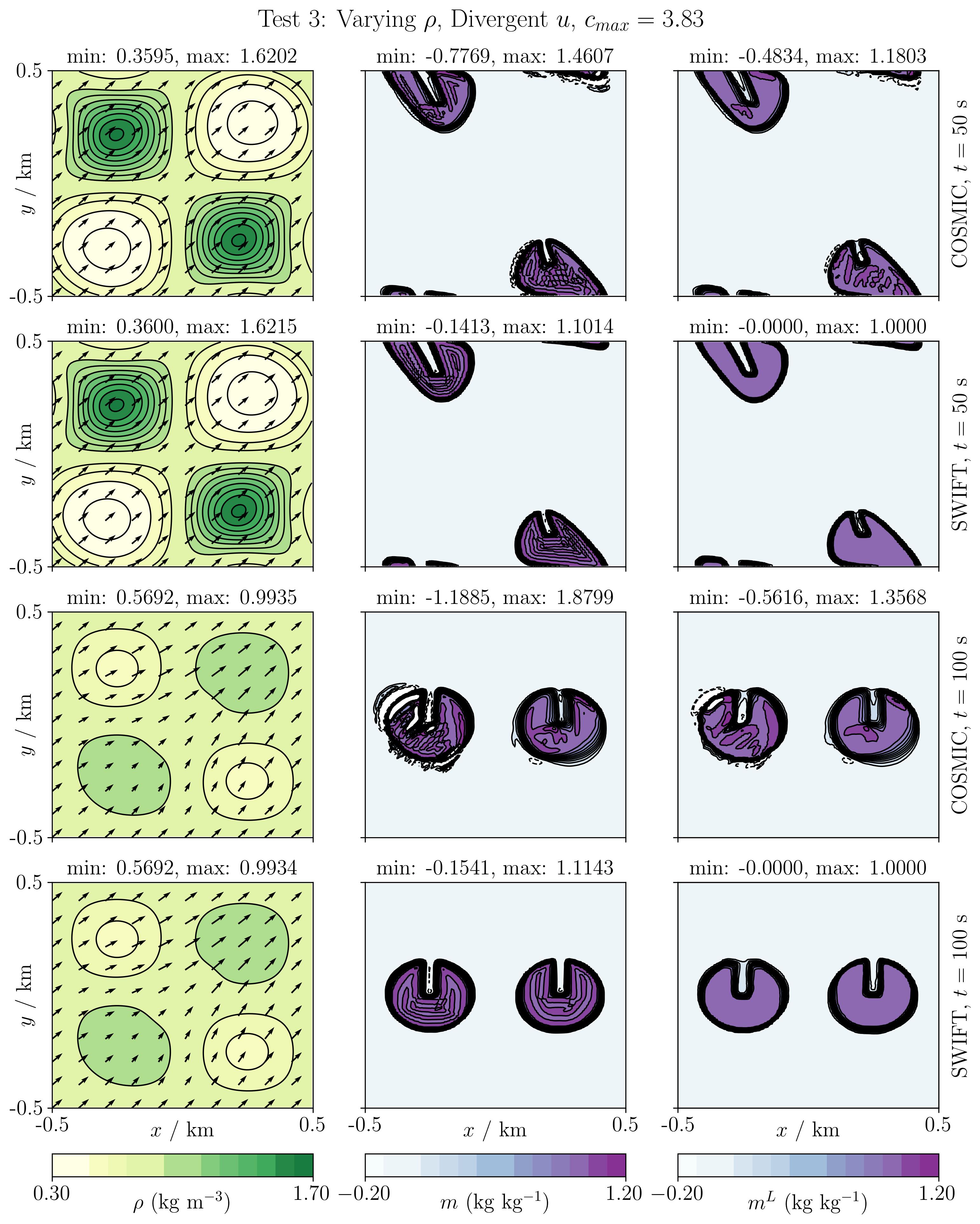}
\caption{The transported fields at the halfway point ($t=50$ s, top two rows) and end of the final time step ($t=100$ s, bottom two rows) of the divergent test of Section \ref{sec:test2}. These fields are taken from the large Courant number case, with $\Delta t=2$ s. Different rows show the COSMIC and SWIFT splittings.
The density fields are contoured at intervals of $0.05$ kg m$^{-3}$, with the $0.8$ kg m$^{-3}$ contour omitted.
The mixing ratio fields have contours every $0.1$ kg kg$^{-1}$, omitting the zero contour.
Arrows on the density plots show the direction and magnitude of the transporting velocity.}\label{fig:test3}
\end{figure}

The non-divergent deformational flow can be modified to make a challenging divergent test case (see Appendix B for the test case setup). The initial density is a smooth field, while the initial tracers are again slotted cylinders. The divergent flow stretches out the tracers, before returning them to their initial conditions.
As with the previous two tests, this was run with two different time steps: $\Delta t=0.2$ s which corresponds to $c_{max}=0.383$, and $\Delta t=2$ s which gives $c_{max}=3.83$. \\
\\
Table \ref{table:test3} contains results from this test.
The table is structured in the same way as Tables \ref{table:test1} and \ref{table:test2}, and demonstrates that with divergent flow the SWIFT splitting still maintains monotonicity properties, while error norms are again comparable with those of the COSMIC splitting.
Fields at $t=50$ s and $t=100$ s from the large Courant number case are plotted in Figure \ref{fig:test3}. \\
 \\
\begin{table}[h!]
\small
\begin{center}
\begin{tabular}{| c | l | c c c | c c c |} \hline
  &  & \multicolumn{3}{c|}{$c_{max}=0.383$} & \multicolumn{3}{c|}{$c_{max}=3.83$} \\
  & Scheme &  Min & Max & $L^2$ error & Min & Max & $L^2$ error   \\ \hline
  \parbox[t]{2mm}{\multirow{6}{*}{\rotatebox[origin=c]{90}{Varying $\rho$}}}
  & COSMIC $\rho$ & 0.597 & 0.997 & 2.24$\times 10^{-3}$ & 0.569 & 0.994 & 2.23$\times 10^{-2}$ \\
  & SWIFT $\rho$ & 0.597 & 0.997 & 2.24$\times 10^{-3}$ & 0.569 & 0.993 & 2.24$\times 10^{-2}$ \\
  & COSMIC $m$ & -0.124 & 1.267 & 2.40$\times 10^{-1}$ & -1.698 & 2.361 & 3.57$\times 10^{-1}$ \\
  & SWIFT $m$ & -0.122 & 1.262 & 2.40$\times 10^{-1}$ & -0.168 & 1.128 & 1.96$\times 10^{-1}$ \\
  & COSMIC $m^L$ & 0.000 & 0.998 & 2.82$\times 10^{-1}$ & -0.597 & 1.409 & 2.70$\times 10^{-1}$ \\
  & SWIFT $m^L$ & 0.000 & 0.998 & 2.80$\times 10^{-1}$ & 0.000 & 1.000 & 2.20$\times 10^{-1}$ \\ \hline
\end{tabular}
\caption{Results for the divergent transport test. The normalised $L^2$ errors, final minimum and maximum values are given for the slotted cylinder initial conditions on the $128\times128$ grid, for Courant numbers $c_{max}=0.383$ and $c_{max}=3.83$.
The mixing ratio $m$ has no limiter applied, while $m^L$ uses the limiter of Appendix A\ref{sec:limiter}.}
\label{table:test3}
\end{center}
\end{table}

\subsection{Convergence Tests} \label{sec:convergence}
\begin{figure}[h!]
\centering
\includegraphics[width=\textwidth]{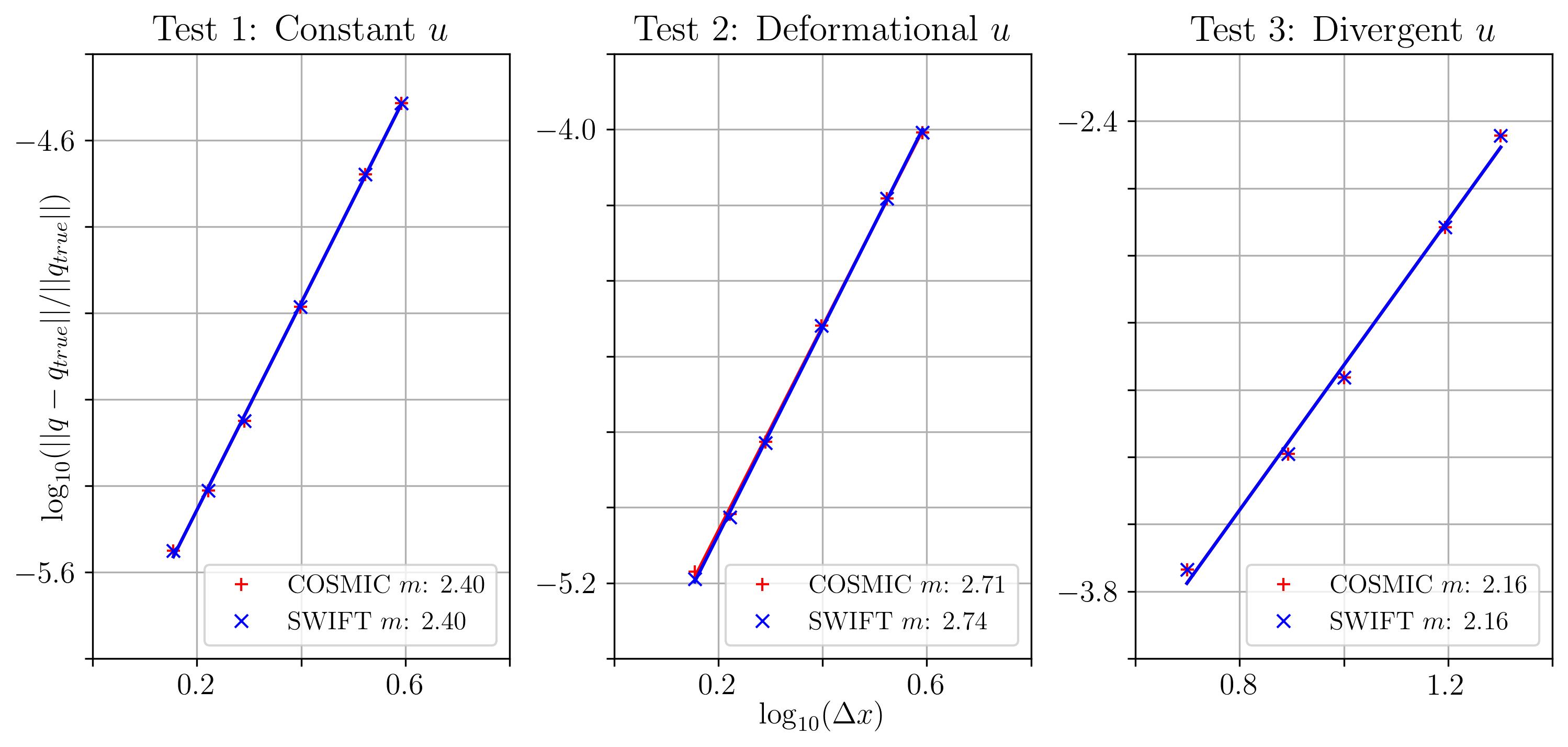}
\caption{The $L^2$ errors in transporting the unlimited tracer mixing ratio, as a function of grid spacing $\Delta x$, whilst keeping $\Delta t$ constant.
Different plots correspond to the different transporting winds described of Section \ref{sec:convergence}.
The legends give the convergence rates, and demonstrate that the error norms for the COSMIC and SWIFT splittings are very similar.
}
\label{fig:convergence_const_dt}
\end{figure}

The test cases of Sections \ref{sec:test1}-\ref{sec:test3} were used to explore the convergence properties of the schemes presented here.
For all convergence tests, the initial tracer fields took the sinusoidal profiles described in Appendix B. \\
\\
We investigated the convergence of the schemes in different regimes.
Using the constant velocity test of Section \ref{sec:test1} and the non-divergent deformation test of \ref{sec:test2}, the relative $L^2$ error norm was measured at different spatial resolutions whilst holding the Courant number constant (by also varying the time step).
This was done for two different maximum Courant numbers ($0.256$ and $2.56$ for test 1, $0.6$ and $6$ for test 2) and on grids with $64^2$, $128^2$, $256^2$ and $512^2$ cells.
Convergence rates for different variables are shown in Table \ref{table:convergence_const_c}.
All convergence rates are determined by fitting a line of best fit between the logarithm of the grid spacing and the relative $L^2$ error.
For the constant velocity (test 1), the convergence rates for the density field and for the tracer field transported relative to a constant density achieve the theoretical third-order values (for the unlimited case).
However when the density is spatially-varying, the tracer field has a second-order convergence rate for both the COSMIC and SWIFT splittings.
We believe that this is because the departure displacement calculation of \eqref{eqn:frac_courant} for the transport of tracers corresponds to an assumption that the density is constant within a cell.
Future work will investigate the impact of using a higher-order subgrid reconstruction of the density to compute the departure displacements.
For the deformational velocity (test 2) with large maximum Courant number, the convergence rates are all approximately second-order. This is because at large Courant numbers with non-constant velocity, the second-order splitting error of the COSMIC and SWIFT splitting dominates. 
\\
\begin{table}[h!]
\small
\begin{center}
\begin{tabular}{| c | l | c  c | c c |} \hline
  &  & Test 1 & & Test 2 &  \\
  &  & $c_{max}=0.256$ & $c_{max}=2.56$ & $c_{max}=0.6$ & $c_{max}=6.0$ \\
  & Scheme &  Convergence rate & Convergence rate & Convergence rate & Convergence rate   \\ \hline
  \parbox[t]{2mm}{\multirow{4}{*}{\rotatebox[origin=c]{90}{Const. $\rho$}}}
  & COSMIC $m$ & 3.01 & 3.01 & 2.37 & 1.98 \\
  & SWIFT $m$ & 3.01 & 3.01 & 2.43 & 1.99 \\
  & COSMIC $m^L$ & 1.87 & 1.78 & 1.69 & 1.98 \\
  & SWIFT  $m^L$ & 1.87 & 1.78 & 1.84 & 1.98 \\ \hline
  \parbox[t]{2mm}{\multirow{6}{*}{\rotatebox[origin=c]{90}{Varying $\rho$}}}
  & COSMIC $\rho$ & 3.01 & 3.01 & 2.37 & 1.98 \\
  & SWIFT $\rho$ & 3.01 & 3.01 & 2.43 & 1.99 \\
  & COSMIC $m$ & 2.00 & 1.99 & 2.06 & 1.99 \\
  & SWIFT $m$ & 2.00 & 1.99 & 2.05 & 1.97\\
  & COSMIC $m^L$ & 1.46 & 1.99 & 1.80 & 1.99 \\
  & SWIFT $m^L$ & 1.38 & 1.99 & 1.84 & 1.96 \\ \hline
\end{tabular}
\caption{The convergence rate of the $L^2$ error with grid spacing, holding the Courant number constant. Errors were measured using the constant velocity test of Section \ref{sec:test1} and the non-divergent deformation test of \ref{sec:test2}, with two different Courant numbers.
The top half of the table gives convergence rates for conservative tracers using a constant density, and the bottom half gives convergence rates for the varying density tests.}
\label{table:convergence_const_c}
\end{center}
\end{table}
\\
Convergence tests were also performed for the tracer with the unlimited PPM scheme, keeping the time step $\Delta t$ constant instead of the Courant number.
This was $\Delta t=0.05$ s for the constant and divergence-free winds, and $\Delta t=0.025$ s for the divergent flow.
Results with each of the different transporting winds from the different test cases are shown in Figure \ref{fig:convergence_const_dt}.
The measured convergence rates at these resolutions lie between 2 and 3, and are comparable with a similar tracer convergence test from \citet{bendall2023solution} that used a Method-of-Lines scheme.
Both COSMIC and SWIFT have very similar error norms at these Courant numbers.

\section{Extensions: Three Dimensions and the Charney-Phillips Grid} \label{sec:3D_CP}
\subsection{Three-Dimensional Scheme} \label{sec:3D}
The extension of the COSMIC splitting to three-dimensions in \cite{leonard1996cosmic} results in multiple flux and advective increments.
To provide a more efficient, yet still high-order splitting we use Strang splitting \citep{strang1968spltting} in the vertical.
For atmospheric applications, breaking the symmetry between the horizontal and vertical dimensions is acceptable given the lack of symmetry between these dimensions in the underlying fluid equations.
However, if COSMIC splitting is used for the horizontal step, the three-dimensional scheme still does not consist of steps with the form \eqref{eqn:mono_condition} and so a monotonic scheme in one-dimensional will not produce a monotonic solution in three-dimensions.
\\
\\
This section will demonstrate how to use the Strang splitting in the vertical with the two-dimensional SWIFT splitting in the horizontal, and thus get the monotonic properties of the one-dimensional scheme for three-dimensional conservative transport.

\subsubsection{Three-Dimensional SWIFT for Density} \label{sec:3D_density}
The first stage is a half step in the vertical.
The corresponding flux calculations use a time step of $\Delta t/2$, which is achieved by using $u^z/2$ in the appropriate flux operator.
\\
\\
The three-dimensional density transport scheme consists of the following steps:
\begin{subequations}
\begin{align}
\rho^z & = \rho^n - \Delta t \nabla_z \cdot \mathcal{F}^z\left(\rho^n, u^z/2, V \right), \\
\rho_I^x & = \rho^z - \Delta t \nabla_x \cdot \mathcal{F}^x\left(\rho^z_A, u^x, \sigma^z V \right), \\
\rho_I^y & = \rho^z - \Delta t \nabla_y \cdot \mathcal{F}^y\left(\rho^z_A, u^y, \sigma^z V \right), \\
\rho^{xy} & = \frac{1}{2}\left[\rho_I^x - \Delta t \nabla_y \cdot \mathcal{F}^y\left(\rho^x_A, u^y, \sigma^xV \right) + \rho_I^y - \Delta t \nabla_x \cdot \mathcal{F}^x\left(\rho_A^y, u^x,\sigma^yV\right)\right], \\
\rho^{n+1} & = \rho^{xy} - \Delta t \nabla_z\cdot\mathcal{F}^z\left(\rho_A^{xy},u^z/2,\sigma^{xy}V \right),
\end{align}
\end{subequations}
where the densities with subscript $A$ are of advective form, and are given by:
\begin{equation}
\rho^z_A = \rho^z / \sigma^z, \quad
\rho^x_A = \rho_I^x / \sigma^x, \quad
\rho^y_A = \rho_I^y / \sigma^y, \quad
\rho^{xy}_A = \rho^{xy} / \sigma^{xy},
\end{equation}
where the transported unity fields are given by
\begin{subequations}
\begin{align}
\sigma^z & = \sigma - \Delta t \nabla_z \cdot \mathcal{F}^z(\sigma,u^z/2,V), \\
\sigma^x & = \sigma^z - \Delta{t} \nabla_x \cdot \mathcal{F}^x(\sigma,u^x,V), \\
\sigma^y & = \sigma^z - \Delta{t} \nabla_y \cdot \mathcal{F}^y(\sigma,u^y,V), \\
\sigma^{xy} & = \sigma^z - \Delta{t} \nabla_x \cdot \mathcal{F}^x(\sigma,u^x,V) - \Delta{t} \nabla_y \cdot \mathcal{F}^y(\sigma,u^y,V).
\end{align}
\end{subequations}
The updated density at the end of the time step can be written as the divergence of fluxes:
\begin{equation}
    \rho^{n+1} = \rho^n - \Delta{t}\nabla_x \cdot F^x  - \Delta{t}\nabla_y \cdot F^y  - \frac{\Delta{t}}{2}\nabla_z \cdot F_1^z - \frac{\Delta{t}}{2}\nabla_z \cdot F_2^z,
\end{equation}
where
\begin{subequations}
\begin{align}
F_1^z & = \mathcal{F}^z(\rho^n,u^z/2,V), \\
F_2^z & = \mathcal{F}^z(\rho^{xy}_A,u^z/2,V), \\
F^x & = \left[ \mathcal{F}^x(\rho_A^z, u^x, V) + \mathcal{F}^x(\rho^y_A, u^x, \sigma^yV) \right]/2, \\
F^y & = \left[ \mathcal{F}^y(\rho_A^z, u^y, V) + \mathcal{F}^y(\rho^x_A, u^y, \sigma^xV)\right]/2.
\end{align}
\end{subequations}
As before, intermediate densities are defined for use in tracer transport:
\begin{equation}
\rho^x = \rho^{n} - \Delta{t} \nabla_x\cdot F^x, \quad
\rho^y = \rho^{n} - \Delta{t} \nabla_y\cdot F^y. 
\end{equation}
If $\rho$ is initially a constant, over a whole time step the increment will be proportional to $\bm{\nabla\cdot u}$, thus satisfying the requirement described by \citet{melvin2024mixed} that density increments should be linear in the divergence.

\subsubsection{Three-Dimensional SWIFT for Tracers}
A three-dimensional conservative scheme for tracers can then be assembled, using Strang splitting between the vertical and horizontal directions and the SWIFT splitting for the horizontal step.
Using variables computed in Section \ref{sec:3D_density}, the scheme can be written as:
\begin{subequations}
\begin{align}
(\rho m)^z & = (\rho m)^n - \Delta t\nabla_z \cdot \mathcal{F}^z\left(m^n, F^z_1, \rho^nV\right), \\
(\rho m)^x & = (\rho m)^z - \Delta t \nabla_x \cdot \mathcal{F}^x\left(m^z, F^x, \rho^zV \right), \\
(\rho m)^y & = (\rho m)^z - \Delta t \nabla_y \cdot \mathcal{F}^y\left(m^z, F^y, \rho^zV \right), \\
(\rho m)^{xy} & = \frac{1}{2}\left[(\rho m)^x - \Delta t \nabla_y \cdot \mathcal{F}^y\left(m^x, F^y, \rho^xV \right) + (\rho m)^y - \Delta t \nabla_x \cdot \mathcal{F}^x\left(m^y, F^x,\rho^yV\right)\right], \\
(\rho m)^{n+1} & = (\rho m)^{xy} - \Delta t \nabla_z\cdot\mathcal{F}^z\left(m^{xy},F_2^z,\rho^{xy}V\right),
\end{align}
\end{subequations}
where $m^x$, $m^y$, $m^z$ and $m_{xy}$ are all found through dividing by the corresponding density fields.
This scheme satisfies all of the desired properties described in Section \ref{sec:intro}: it is conservative for the tracer whilst preserving a constant field, and the multi-dimensional scheme inherits any monotonicity properties of the underlying one-dimensional scheme.
These properties are not compromised for flows with large Courant number.

\subsection{The Charney-Phillips Staggering} \label{sec:charney-phillips}
\begin{figure}[h!] 
\centering
\includegraphics[width=0.8\textwidth]{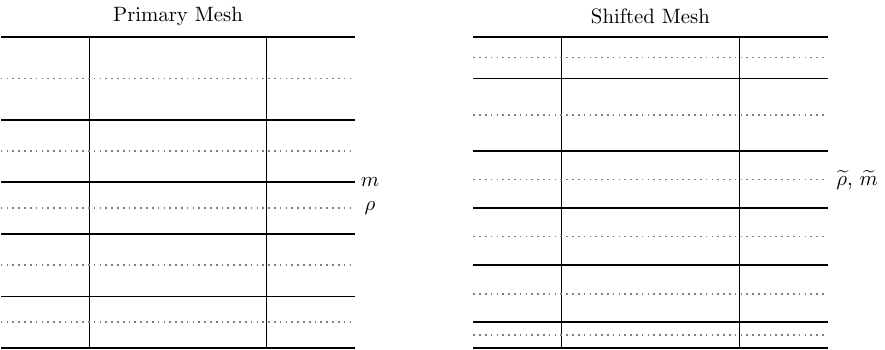}
\caption{A vertical cross-section of an example vertically-shifted mesh, used to transport tracers on the Charney-Phillips grid.
Solid black lines represent the surfaces between model layers, while dotted grey lines showing the vertical centres of the layers.
With the Charney-Phillips staggering on the primary mesh, the mixing ratio $m$ is located on the surfaces between model layers on the primary mesh, while the dry density is located at the centres of layers.
On the vertically-shifted mesh, these quantities are described at cell centres.
The surfaces between model layers on the vertically-shifted mesh are at the midpoints of layers on the primary mesh, so that the layers are shifted relative to those on the primary mesh.
The vertically-shifted mesh therefore has one more level than the primary mesh.
} \label{fig:shifted}
\end{figure}

The previous sections have focused on a tracer transport scheme for the case that the mixing ratio $m$ is co-located with the dry density $\rho$.
However, for dynamical cores with a Charney-Phillips staggering, mixing ratio fields may be vertically-staggered with $\rho$.
For instance these variables could be moisture species, which are instead co-located with the entropy-type variable (such as potential temperature) in order to accurately capture the latent heat exchanges associated with phase changes.
It is still desirable for the transport of these staggered tracers to possess the same properties as those discussed in Section \ref{sec:intro}: local conservation, stability independent of Courant number, preservation of a constant and options for monotonicity or positivity.
This section summarises how these properties can still be attained with the FFSL scheme described in the previous sections, by using the framework of \citet{bendall2023solution}.
\\
\\
The problem of achieving consistent and conservative transport of moisture species on the Charney-Phillips grid was addressed by \citet{bendall2023solution}, in the context of a Method-of-Lines transport scheme.
The approach used there was to express the corresponding density of the tracer species on a mesh that is \textit{vertically-shifted} from the primary mesh used by the model.
This shifted mesh was also introduced by \citet{thuburn2022numerical}.
On the shifted mesh, the surfaces between model layers are at the midpoints of the model layers of the primary mesh.
The top and bottom surfaces of the shifted mesh coincide with those of the primary mesh, so that the lowest and highest layers are of half the depth of other layers.
Then the shifted mesh has one extra layer compared with the primary mesh, which ensures that variables expressed at cell centres on the shifted mesh have the same number of degrees of freedom as those described on the surfaces between model layers on the primary mesh.
This is illustrated in Figure \ref{fig:shifted}.
\\
\\
In the framework of \citet{bendall2023solution}, the conservative transport of tracers takes place on this shifted mesh.
This involves mapping $\rho$ and $m$ from the primary mesh to the centres of cells on the shifted mesh, and also mapping the density flux $\bm{F}$ from facets on the primary mesh to facets on the shifted mesh.
As in \citet{bendall2023solution}, in this work quantities on the shifted mesh are denoted by a tilde $\widetilde{\cdot}$ adornment. \\
\\
The tracer mass in the $i$-th cell on the shifted mesh is given by
\begin{equation}
\widetilde{\rho}_i\widetilde{m}_i \widetilde{V}_i,
\end{equation}
where $\widetilde{V}_i$ is the volume of the shifted mesh cell.
Conservative tracer transport therefore conserves the mass that is the global sum of this quantity, whilst solving a transport equation like \eqref{eqn:tracer_conservative} ensures that this mass is locally-conserved.
There are three extra ingredients to obtain the consistency property.
Firstly, as in Sections \ref{sec:consistent_ffsl_1d}  and \ref{sec:swift_tracer} of this paper, the flux that was used in the transport of $\rho$ must be re-used in the transport of the tracer (but on the shifted mesh).
Secondly, the operators to remap fields from the primary mesh to the shifted mesh should satisfy some particular conditions given in \citet{bendall2023solution}, which guarantee that the processes of transporting a density and converting it to the shifted mesh commute with one another.
In this work, we therefore use the remapping operators presented in \citet{bendall2023solution}, with one amendment: to avoid generating new extrema in the top and bottom half-layers, the transformation between $m$ on the primary and shifted meshes is simply the identity matrix, rather than performing interpolation/extrapolation at the domain boundaries.
Finally, the operator to evaluate the flux should satisfy:
\begin{equation}
\mathcal{F}(K,F,\rho V) = KF,
\end{equation}
when the tracer mixing ratio is a constant $K$.
As described in Section \ref{sec:consistent_ffsl_1d}, the tracer fluxes in this work do satisfy this property, which ensures that the transport of $m$ will be consistent despite the Charney-Phillips staggering.

\subsection{Three-Dimensional Results} \label{sec:results3d}
To demonstrate the performance of SWIFT in three-dimensions we use the three-dimensional deformational test of \cite{skamarock2006limiters} with an added background velocity in the horizontal (see Appendix B). The domain is doubly periodic in the horizontal, with zero flux boundaries at the model top and bottom. The tests are performed on the $64\times64\times64$ grid, with step function initial tracers and a smooth density field. We consider two tracers, $m_c$ co-located with $\rho$ and $m_s$ vertically-staggered with $\rho$.\\
\begin{table} [h!]
\small
\begin{center}
\begin{tabular}{| l | c c c | c c c |} \hline
  & \multicolumn{3}{c|}{$c_{max}=0.48$} & \multicolumn{3}{c|}{$c_{max}=4.8$} \\
Scheme & Min & Max & $L^2$ error & Min & Max & $L^2$ error \\ \hline
COSMIC $\rho$ & 0.504 & 0.996 & 8.19$\times 10^{-5}$ & 0.504 &  0.996 & 9.28$\times 10^{-4}$ \\ 
SWIFT $\rho$ & 0.509 & 0.996 & 8.18$\times 10^{-5}$ & 0.504 & 0.996 & 9.47$\times 10^{-4}$ \\
COSMIC $m_c$ & -0.298 & 1.216 & 1.75$\times 10^{-1}$ & -0.296 & 1.352 & 1.61$\times 10^{-1}$ \\ 
SWIFT $m_c$ & -0.299 & 1.214 & 1.75$\times 10^{-1}$ & -0.136 & 1.140 & 1.54$\times 10^{-1}$ \\ 
COSMIC $m^L_c$ & 0.000 & 1.000 & 2.28$\times 10^{-1}$ & -7.6$\times 10^{-3}$ & 1.100 & 1.91$\times 10^{-1}$  \\ 
SWIFT $m^L_c$ & 0.000 & 1.000 & 2.27$\times 10^{-1}$ & 0.000 & 1.000 & 1.90$\times 10^{-1}$ \\  
COSMIC $m_s$ & -0.304 & 1.215 & 1.70$\times 10^{-1}$ & -0.296 & 1.352 & 1.49$\times 10^{-1}$ \\ 
SWIFT $m_s$ & -0.305 & 1.214 & 1.69$\times 10^{-1}$ & -0.134 & 1.144 & 1.41$\times 10^{-1}$ \\ 
COSMIC $m^L_s$ & 0.000 & 1.000 & 2.16$\times 10^{-1}$ & -5.9$\times 10^{-3}$ &  1.100 & 1.77$\times 10^{-1}$ \\ 
SWIFT $m^L_s$ &  0.000 & 1.000 & 2.16$\times 10^{-1}$ & 0.000 &  1.000 & 1.77$\times 10^{-1}$ \\ 
\hline
\end{tabular}
\caption{Final minima, maxima, and $L^2$ error norms for the three-dimensional test with slotted cylinder initial tracers and smooth varying density. Results are presented for two different maximum Courant numbers, $4.8$ and $0.96$. Subscript $c$ indicates that the tracer $m$ is co-located with $\rho$, subscript $s$ indicates the tracer is staggered with respect to $\rho$ in the vertical. Results are from the $64\times64\times64$ grid. Values are very close between COSMIC and SWIFT for the small Courant number case (left-hand side). For the large Courant number test (right-hand side), with the limiter the SWIFT splitting is monotonic for both the co-located and staggered tracers.}
\label{table:test3d}
\end{center}
\end{table}
\begin{figure}[ht!]
\centering
\includegraphics[width=0.9\textwidth]{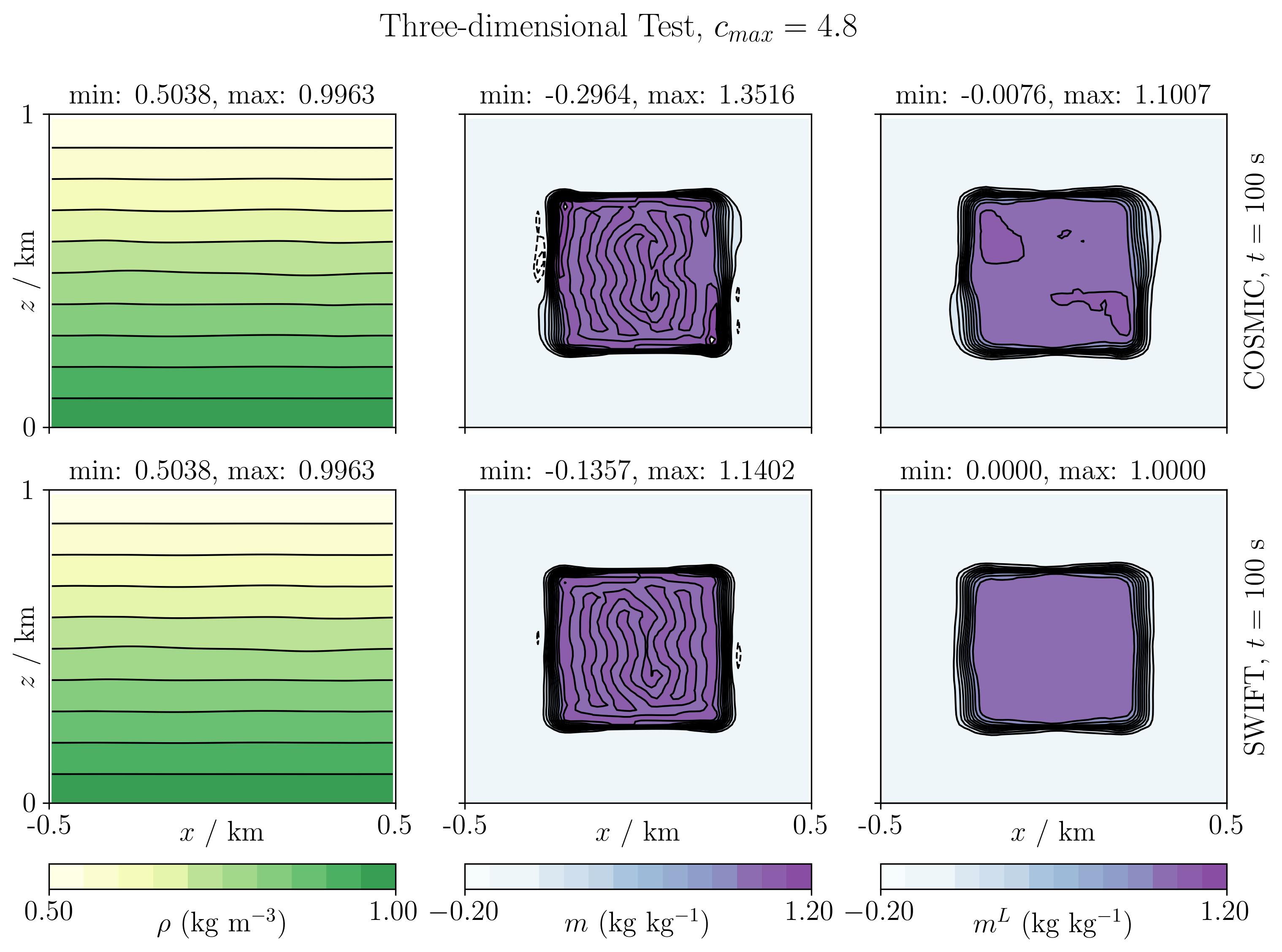}
\caption{Slices along the $x$-$z$ plane of transported fields from the three-dimensional transport test of Section \ref{sec:results3d}. These fields are taken from the large Courant number case, with $\Delta t=2.5$ s. The top row shows fields transported with the COSMIC splitting, while the bottom row contains fields transported with the SWIFT splitting.
Mixing ratios are vertically-staggered from the density field, demonstrating that the splitting performs as expected on the Charney-Phillips grid.
The density fields are contoured at intervals of $0.05$ kg m$^{-3}$.
The mixing ratio fields have contours every $0.1$ kg kg$^{-1}$, omitting the zero contour.}
\label{fig:test3d}
\end{figure} \\
Table \ref{table:test3d} shows the results for the three-dimensional test for two different maximum Courant numbers: $c_{max}=4.8$ using $\Delta t=2.5$ s, and $c_{max}=0.48$ using $\Delta t=0.25$ s. Note that this test has different magnitudes of velocity in each direction (due to the deformational coefficient and the background flow) and so the maximum Courant number is in the $x$ direction. For the small Courant number test the SWIFT and COSMIC schemes give near identical results. Using the strict limiter gives monotonicity for both SWIFT and COSMIC in this case.
\\
\\
For the large Courant number case the error norms are similar for both SWIFT and COSMIC. However, the results show that COSMIC is not monotonic even when the strict limiter is used for PPM. The SWIFT scheme with the strict limiter is monotonic, for tracers co-located or staggered with $\rho$ in the vertical.
Vertical slices of the final fields are shown in Figure \ref{fig:test3d}.

\section{Summary}

It is desirable that transport schemes for tracers in atmospheric dynamical cores are conservative, consistent with the transport of the dry density, long time-step-permitting, and can be monotonic or positive-definite.
Flux-Form Semi-Lagrangian (FFSL) schemes can be designed to have these properties in one-dimension, and splitting techniques can be used to solve multi-dimensional transport problems.
The splitting of \cite{lin1996ffsl,leonard1996cosmic}, known as the COSMIC splitting, allows for consistent, long time-step-permitting two-dimensional transport, but the monotonicity or positivity properties of the one-dimensional scheme are not replicated for the two-dimensional scheme. \\
\\
Here we have presented the SWIFT scheme, which builds upon existing FFSL schemes, the COSMIC splitting, and the work of \cite{skamarock2006limiters}, to produce a two-dimensional splitting that retains the monotonicity or positivity properties of the one-dimensional schemes, even for large Courant number flow. The SWIFT splitting can be used for consistent tracer transport, the transport of a density, or even for the advective transport of a tracer. \\
\\
The use of the two-dimensional SWIFT splitting in the horizontal can be combined with vertical transport through a Strang splitting.
This provides a three-dimensional SWIFT scheme which still retains all of the desirable properties discussed above.
The three-dimensional SWIFT scheme can also be applied to tracer transport on a staggered vertical grid, such as the Charney-Philips, in which the dry density and the tracer mixing ratio are not co-located. \\
\\
The SWIFT scheme has been compared with the COSMIC splitting for a number of numerical tests.
The Piecewise Parabolic Method (PPM, \cite{colella1984ppm}) is used as the one-dimensional scheme in each case, with the unlimited and monotonic limited versions being tested.
For the two-dimensional tests, SWIFT splitting and COSMIC splitting are shown to be identical when using a constant velocity and constant density and unlimited PPM. In the case of the varying wind or varying density test cases, the results show that the SWIFT splitting generally has a smaller error norm than the corresponding COSMIC splitting, and confirms that the SWIFT splitting maintains any monotonicity or positivity of limiters used for PPM, whilst COSMIC splitting does not.
The tests also demonstrate that this holds for large Courant number flows. The three-dimensional testing again shows that SWIFT splitting retains these properties when Strang splitting is used in the vertical, even when the tracer is staggered relative to the dry density. \\
\\
The benefits of the scheme presented here have seen it become a candidate for the transport scheme used by the Met Office's new GungHo dynamical core \citep{melvin2024mixed}.
In future work, we will describe how this scheme can be applied to spherical geometries, and in particular cubed-sphere meshes.
This scheme can also be combined with the work of \citet{thuburn2022numerical} to achieve entropy conservation on the Charney-Phillips grid, or with the approach of \citet{brown2024physics} to transport tracers which are expressed on a coarser mesh than that used by the dynamical core.

\section*{Acknowledgements}
The authors would like to acknowledge the frequent discussions with Mohamed Zerroukat, Ben Shipway, Thomas Melvin and Nigel Wood which contributed to shaping this work, and for their feedback on early drafts of the manuscript. The authors are grateful to Hilary Weller whose suggestions helped improve the manuscript, and also to Lucas Harris and the two anonymous referees for their reviews which led to an improved text.

\section*{Data Availability}
The results contained in this work were produced with the LFRic-Atmosphere model.
The LFRic-Atmosphere source code and configuration files are freely available from the Met Office Science Repository Service \\
(https://code.metoffice.gov.uk) upon registration and completion of a software licence.
Plotting scripts are stored in a  public Github repository (https://github.com/tommbendall/bendall\_kent\_swift\_2024\_plots).

\appendix

\section{The Piecewise Parabolic Method} \label{sec:ppm_appendix}

This appendix describes the Piecewise Parabolic Method (PPM, \citet{colella1984ppm}) used for the fractional flux in the testing in this paper.

\subsection{PPM Coefficients} \label{sec:ppm_recon}

The unlimited PPM scheme \citep{colella1984ppm} computes the quadratic reconstruction such that: 1) the integral of the reconstruction equals the integral of the field within the departure cell; and 2) the reconstruction at the centres of cell facets equals a high-order interpolation of the field at those facets. \\
\\
Let $q_0$ and $q_1$ be the interpolation of a scalar field $q$ to cell $i$ facets. Setting $\mathcal{Q}(0)=q_0$ and $\mathcal{Q}(1) = q_1$ in \eqref{eqn:quad_coef} (i.e. at $\xi=0$ and $\xi=1$ respectively), and setting the integral of the subgrid reconstruction equal to $q_i$, gives the PPM coefficients
\begin{subequations}
\begin{align}
    a_0 &= q_0, \\
    a_1 &= -4 q_0 - 2 q_1 + 6 q_i, \\
    a_2 &= 3 q_0 + 3 q_1 - 6 q_i.
\end{align}
\end{subequations}
The edge values are computed using a high-order interpolation, usually fourth-order \citep{colella1984ppm} or sixth-order \citep{colella2008}.
For a uniform grid the fourth-order facet interpolation is
\begin{equation}
    q_{i-1/2} = q_0 = \left( -q_{i-2} + 7 q_{i-1} + 7 q_{i} - q_{i+1}\right)/12.
\end{equation}
As the coefficients only depend on the field in the departure cell and the interpolation to cell facets, for non-uniform grids the coefficients are the same as above, with only the interpolation being recalculated.  

\subsection{Monotonic Limiter} \label{sec:limiter}
A monotonic limiter can be applied to the quadratic subgrid reconstruction to prevent unphysical amplification of extrema.
For a tracer mixing ratio, this should ensure that the field remains bounded by its initial extrema.
Such a monotonic limiter will also enforce positivity in one dimension. 
Here we consider one example of a limiter, which works by ensuring that the reconstruction does not exceed the bounds set by the neighbouring grid cells.  \\
\\
First, the stationary point of the quadratic is tested. Let $\tau$ be the stationary point, which is found from
\begin{equation}
a_1 + 2 a_2 \tau  = 0 \Rightarrow \tau = - \frac{a_1}{2 a_2}. 
\end{equation}
If $\tau(1-\tau) > 0$ then we revert to a constant reconstruction
\begin{subequations}
\begin{align}
    a_0 &= q_i, \\
    a_1 &= 0, \\
    a_2 &= 0.
\end{align}
\end{subequations}
For PPM, the facet values are computed by a high-order interpolation. If a facet value exceeds the neighbouring field values, we force the facet value to lie between them, e.g. for $q_1$ if $(q_1 - q_i)( q_{i+1} - q_1) < 0$, we set the facet value $q_1 = \min[ \max(q_{i+1},q_i), \max(q_1, \min[q_{i+1},q_i]) ]$. Note that the limiter presented here is stricter than that given in \citep{colella1984ppm}, resulting in a more diffusive solution. 

\subsection{PPM Reconstruction} \label{sec:consistent_ppm}

Integrating the PPM coefficients and rearranging gives the PPM facet reconstruction used for the fractional fluxes in the form of \eqref{eqn:consistent_flux}.
As in Appendix A\ref{sec:ppm_recon}, $q_0$ and $q_1$ are interpolations of $q$ to the centres of cell facets.
Let cell $i$ be the upwind cell from the desired facet, then
\begin{equation}
\mathcal{R}(q,c) =
\begin{cases}
(1-2c+c^2)q_1 + (3c-2c^2)q_{i} + (-c+c^2)q_0, & c \geq 0, \\
(c+c^2)q_1 + (-3c-2c^2)q_{i} + (1+2c+c^2)q_0 & c < 0,
\end{cases}
\end{equation}
where for brevity, $c=c^{frac}_{i+1/2}$. \\
\\
The PPM limiter checks for extrema in the reconstruction (note that the interpolated edge values must also be monotonic). The limiter sets
\begin{equation}
\tau = (2q_0+q_1-3q_{i})/(3q_0+3q_1-6q_{i}),
\end{equation}
then if $\tau (1-\tau) > 0$ the reconstruction reverts to first-order by setting the reconstruction to be the field value in the upwind cell
\begin{equation}
\mathcal{R}(q,c) = q_i.
\end{equation}
With this reconstruction, if the field is some constant $K$, then $q_0=q_i=q_1=K$, and $\mathcal{R}(q,c)=K$.

\section{Test Cases} \label{sec:test_appendix}

This section describes the test cases used for the numerical testing in this article. 

\subsection{Two-Dimensional Tests}

The two-dimensional tests are set up on a doubly-periodic planar domain, with $500 \leq x,y \leq 500$ m, so that the domain's lengths were $L_x=1000$ m and $L_y=1000$ m.
Tests are run to a time of 100 s.
\\
\\
The tests either use a constant density, with $\rho=1$ kg m$^{-3}$, or a smooth varying density given by
\begin{equation}
    \rho = \rho_{\mathrm{ref}} + \rho_{\mathrm{mag}} \sin(2 \pi x/L_x) \sin(2 \pi y/L_y),
\end{equation}
where $\rho_{\mathrm{ref}}=0.8$ kg m$^{-3}$ and $\rho_{\mathrm{mag}}=0.2$ kg m$^{-3}$. \\ 
\\
The tracers are either set up as sine waves or slotted cylinders. The sine wave is given by
\begin{equation}
    m = m_{\mathrm{ref}} + m_{\mathrm{mag}} \sin(2 \pi x/L_x) \sin(2 \pi y/L_y),
\end{equation} 
where $m_{\mathrm{ref}}=m_{\mathrm{mag}}=0.5$ kg kg$^{-1}$. \\
\\
For the slotted cylinders we set the centre of the tracers to be $x_{c1} = -250$ m and $x_{c2}=250$ m, and $y_c=0$ m.
The slotted cylinder initial condition requires the distance from the centre to every point in the domain
\begin{equation}
    r_i = \sqrt{(x-x_{ci})^2+(y-y_c)^2}, \hspace*{1cm} i=1,2.
\end{equation}
The slotted cylinders are set up as
\begin{equation}
    m = 
    \begin{cases}
        m_c, & \mathrm{if} \ r_i < r_c, \\
        0, & \mathrm{if} \ y>y_c \ \mathrm{and} \ |x-x_{ci}|<l_c, \\
        0, & \mathrm{otherwise}.
    \end{cases}
\end{equation}
with $m_c=1$ kg kg$^{-1}$, $r_c=160$ m and $l_c=25$ m.
For the constant wind test the velocity is set as $u^x=10$ m s$^{-1}$ and $u^y=10$ m s$^{-1}$. For the non-divergent deformational flow we follow \cite{skamarock2006limiters} and \citet{kent2020positive} but scale up to our domain size to give the velocities as
\begin{subequations}
\begin{align}
u^x & =  u_0 \sin^2(\pi x'/L_x)\sin(2 \pi y'/L_y)\cos(\pi t/T)+u_0,\\ 
u^y & = -u_0 \sin^2(\pi y'/L_y)\sin(2\pi x'/L_x)\cos(\pi t/T)+u_0,
\end{align}
\end{subequations}
where $u_0=10$ m s$^{-1}$, $T=100$ s, and the moving coordinates are given by $x' = (x+L_x/2) - u_0 t$ and $y' = (y+L_y/2) - u_0 t$. The divergent deformational flow is given by
\begin{subequations}
\begin{align}
u^x & =  \tfrac{1}{2}u_0 \sin^2(\pi x'/L_x)\sin(2 \pi y'/L_y)\cos(\pi t/T)+u_0,\\ 
u^y & = \tfrac{1}{2}u_0 \sin^2(\pi y'/L_y)\sin(2 \pi x'/L_x)\cos(\pi t/T)+u_0.
\end{align}
\end{subequations}

\subsection{Three-Dimensional Tests}

The three-dimensional tests are set up on planar domain, periodic in $x$ and $y$, with $-500 \leq x,y \leq 500$ m and $0 \leq z \leq L_z$, so that the domain lengths are $L_x=L_y=L_z=1000$ m. \\
\\
The density is a linear function of height
\begin{equation}
    \rho = \rho_{\mathrm{ref}} + \rho_{\mathrm{mag}}(1-z/L_z),
\end{equation}
where $ \rho_{\mathrm{ref}}=\rho_{\mathrm{mag}} = 0.5$ kg m$^{-3}$.  \\
\\
The tracer is initialised as a step function with $m_c=1$ kg kg$^{-1}$ so that
\begin{equation}
    m = 
    \begin{cases}
        m_c, & \mathrm{if} \ |x| < L_x/4 \ \mathrm{and} \ |z-L_z/2|< 3L_z/10, \\
        0, & \mathrm{otherwise}.
    \end{cases}
\end{equation}
For the three dimensional test the velocity is taken from \cite{skamarock2006limiters}. Scaling up to our domain size gives
\begin{subequations}
\begin{align}
u^x & =  2 u_0 \sin^2(\pi x'/L_x)\sin(2 \pi y'/L_y)\sin(2 \pi z/L_z)\cos(\pi t/T)+u_0,\\ 
u^y & = -u_0 \sin^2(\pi y'/L_y)\sin(2 \pi x/L_x)\sin(2 \pi z/L_z)\cos(\pi t/T)+u_0,\\ 
u^z & = -u_0 \sin^2(\pi z/L_z)\sin(2 \pi x/L_x)\sin(2 \pi y/L_y)\cos(\pi t/T),
\end{align}
\end{subequations}
where $u_0=10$ m s$^{-1}$ and $x'$ and $y'$ are defined as for the two-dimensional deformational case.

\bibliography{main.bib}

\end{document}